\documentclass{article}
\usepackage[ansinew]{inputenc}
\usepackage[english]{babel}
\usepackage{graphicx}
\usepackage{amsmath,amssymb,amsthm,mathrsfs,amsfonts}

\newtheorem {thm}{Theorem}[section]
\newtheorem {defn}{Definition}[section]
\newtheorem {prop}[thm]{Proposition}

\newtheorem {lemma}[thm]{Lemma}

\newtheorem {cor}[thm]{Corollary}
\newtheorem{conj}{Conjecture}
\newcommand{\Ir}{\boldsymbol{\mathcal{I}_r}}

\newcommand{\SP}{{\operatorname{sp}}}
\newcommand{\TR}{{\operatorname{tr}}}
\newcommand{\SPEC}{{\operatorname{spec}}}
\newcommand{\bbmatrix}[1]{\left[ \begin{array}{cccccccccccccccccc} #1 \end{array} \right]}

\title{On a theory of vessels and the inverse scattering}
\author{Andrey Melnikov\\Drexel university, Philadelphia, USA}
\begin{document}
\maketitle

\abstract{In this paper we present a theory of vessels and its application to the classical inverse scattering
of the Sturm-Liouville differential equation.
The classical inverse scattering theory, including all its ingredients: Jost solutions, the Gelfand-Levitan
equation, the tau function, corresponds to regular vessels, defined by bounded operators.
A contribution of this work is the construction of models of vessels corresponding to unbounded operators, which
is a first step for the inverse scattering for a wider class of potentials. 

A detailed research of Jost solutions and the corresponding vessel is presented for the 
unbounded Sturm-Liouville case.
Models of vessels on curves, corresponding to unbounded operators are presented as a tool to study Linear Differential
equations of finite order with a spectral parameter and as examples, we show how the family of Non Linear
Schr\" odinger equations and Canonical Systems arise.
}

\tableofcontents

\section{Introduction}
The Sturm Liouville differential equation \cite{bib:Liouville, bib:Sturm} is one of the oldest differential equations,
studied by mathematicians. It is defined as a linear differential equation of second order\begin{equation}\label{eq:SL}
-\frac{d^2}{dx^2} y(x) + q(x) y(x) = s^2 y(x),
\end{equation}
where $\lambda\in\mathbb C$ is called the \textit{spectral parameter} and the coefficient $q(x)$
is called the \textit{potential}.
It is the simplest linear equation, for which one can not usually find closed-form solutions.

This equation was and probably is one of the most influential with mathematical analysis,
because many techniques were developed in order to solve it. 
For example, it was studied 
\begin{enumerate}
	\item by C. Sturm \cite{bib:Sturm}, and R. Liouville \cite{bib:Liouville} in connection with the dynamics, the heat
		equation,
	\item using the monodromy preserving deformation problem of Linear Differential Equations (LDEs)
		by L. Schlesinger \cite{bib:Schlez}, R. Fuchs \cite{bib:Fuchs} and Garnier \cite{bib:Garnier},
	\item using the Scattering theory by Lax--Phillips \cite{bib:LaxPhil}, and Gelfand-Levitan \cite{bib:GL},
	\item using Riemannian transformations by A. Povzner \cite{bib:Povzner} and V.A. Marchenko 
		\cite{bib:Marchenko,bib:MarchenkoSL}.
\end{enumerate}
Also M. G. Krein \cite{bib:Krein1} and many other famous mathematicians gave fundamental contributions to
this equation. Actually, the list of the contributors and techniques can easily fill few pages.

The third theory in this list, the Scattering theory, studies asymptotic behavior of solutions of the equation 
\eqref{eq:SL} and compares them to the trivial ones, corresponding to the zero potential.
Notice that solving SL equation \eqref{eq:SL} for 
$q(x)=0$ one obtains that the solutions are linear combinations of the exponents $e^{isx}, e^{-isx}$.
For the potential $q(x)$, which is locally integrable and satisfies the condition \cite{bib:FadeevInv}
\begin{equation} \label{eq:FadAssum}
\int\limits_0^\infty x |q(x)| dx = C < \infty.
\end{equation}
one can define \textit{Jost} solutions, which behave asymptotically (when $x \rightarrow\infty$) as the 
trivial ones with a certain phase, depending on the
spectral parameter $s$. Following L.D. Faddeyev \cite{bib:FadeevInv} 
\textit{"the fundamental problem arising in the quantum theory of scattering
is the solution of
\[ L \psi(x,k) = -(\dfrac{\partial^2}{\partial x^2} \psi(x,k)) + q(x) \psi(x,k) = k^2 \psi(x,k)
\]
satisfying the condition $\psi(0,k) = 0$, behaves asymptotically like 
$\psi(x,k)\approx C(k) \sin(kx-\eta(k))$ provided the potential $q(x)$ decreases sufficiently fast as x tends to
infinity; to what extent does the assignment of $\eta(k)$ determine the potential $q(x)$ and how these functions
are related"}.

In this paper, we will present a theory, which generalize the idea of the inverse scattering, i.e. which finds a correspondence between potentials and some (matrix-valued) functions of a complex variable in a slightly different setting, and which coincides with the classical
inverse scattering in a "regular" case. The benefit is that we can unify all the approaches
and apply this theory to study for example NLS equations (section \ref{sec:NLS}) and Canonical systems (section \ref{sec:CanSys})  beyond the classical results. It is important to notice that
it is a separate project by itself and the present work is a background for this future work.

Under slightly different assumption on the potential \cite{bib:FaddeyevII}
$\int\limits_{-\infty}^\infty (1+|x|) |q(x)| dx = C < \infty$ one can
construct the Jost solutions $f_1(s,x)$ and $f_2(s,x)$ of \eqref{eq:SL} such that
\[ \lim\limits_{x\rightarrow\infty} \dfrac{f_1(s,x)}{e^{isx}} = 
\lim\limits_{x\rightarrow-\infty} \dfrac{f_2(s,x)}{e^{-isx}} = 1
\]
and create the following matrix of $\lambda=\sqrt{is^2}$, where we choose $\Im s\geq 0$
and $\beta=\dfrac{1}{2}\int q$:
\[ S(\lambda,x) = \bbmatrix{\dfrac{f_1 + f_2}{2}& \dfrac{f_1 - f_2}{2s}\\
\dfrac{f_1' + f_2'}{2i}-\beta\dfrac{f_1 + f_2}{2i} & \dfrac{f_1' - f_2'}{2is}-\beta\dfrac{f_1 - f_2}{2is}} \bbmatrix{\cos(sx)&\dfrac{i\sin(sx)}{s}\\is\sin(sx)&\cos(sx)}^{-1},
\]
which has 4 defining properties
\begin{enumerate}
	\item behaves like $I$ for $\lambda$ approaching infinity, and has a jump along a 
		finite cut $\Gamma$ on the imaginary positive axis,
	\item twice differentiable with respect to $x$,
	\item symmetric (with respect to Pauli matrix $\sigma_1=\bbmatrix{0&1\\1&0}$) see
		\eqref{eq:Symmetry},
	\item maps solutions of the trivial SL equation to solutions of the complicated one,
		if one concentrates on the first row.
\end{enumerate}
These properties are similar to the properties of the scattering matrix, appearing
in \cite[section 5]{bib:FaddeyevII}. It turns out that one can realize this matrix,
using \cite{bib:bgr} in the form (for each $x$)
\[ S(\lambda,x) = I - B^*(\mu, x) \mathbb X^{-1}(x) (\lambda - A)^{-1} B(\mu, x) \sigma_1,
\]
where there arises an auxiliary separable Hilbert space $\mathcal H$ and operators
$A, \mathbb X(x): \mathcal H \rightarrow \mathcal H$,
$B(x): \mathbb C^2 \rightarrow \mathcal H$, satisfying certain relations (see 
Definition \ref{def:KV}). One can also represent such a function in the form
\[ S(\lambda,x) = I + \int\limits_{it\in\Gamma} \dfrac{S(\mu(t), x)}{\lambda - it} dt
\]
Since $S(\lambda,x) - I$ satisfies conditions of the limiting values on an axis theorem, i.e.
$S(\lambda,x)$ is represented as a Poisson integral of its limiting values on the axis.
These two ideas brought us to a far reaching generalization. It turns out that one can
generalize the construction of such a matrix $S(\lambda,x)$ not only for SL equation but also
to a wider class of differential equations (this is done in section \ref{sec:DefVel}).

The background for this research is the work \cite{bib:SchurVessels} (which was announced in \cite{bib:amv})
and a realization theory of matrix-valued $p\times p$
functions of a complex variable $\lambda$, analytic and invertible (hence identity) at infinity, 
and $J$-contractive ($J=J^*=J^{-1}$) \cite{bib:BL, bib:KreinReal}. On the one hand, they have
a so called \textit{realization theorem} (based on Theorem \ref{thm:RealizeSinK})
\[ S(\lambda) = I -  B^* \mathbb X^{-1} (\lambda I - A)^{-1} B J,
\]
where $A, \mathbb X$ are selfadjoint bounded operators, acting on an auxiliary Krein space $\mathcal H$ and 
$B:\mathbb C^p \rightarrow \mathcal H$ is also bounded. On the other hand one can apply a "vessel construction" (see Section \ref{sec:StandConstr}) 
and to obtain a vessel, whose transfer function depends additionally on a real variable $x$ and is of the form
\[ S(\lambda,x) = I - B^*(x) \mathbb X^{-1}(x)(\lambda I - A)^{-1} B(x) \sigma_1(x).
\]
It holds that for
$x=x_0$ the function $S(\lambda,x)$ coincides with the above realization for $S(\lambda)$.

Starting from $S(\lambda)$ realized with bounded operators, from the properties of the vessel
it follows that the class of the transfer functions
\[ S(\lambda,x) \in \Ir(\sigma_1(x), \sigma_2(x), \gamma(x), \gamma_*(x), \mathrm I), \]
of so called regular vessels, constructed in this manner, consists of
functions which are $\sigma_1(x)$ symmetric, identity around $\lambda=\infty$ and
map solutions of the \textit{input} Linear Differential Equation (LDE) \eqref{eq:InCC} with spectral parameter
$\lambda$
\[ -\sigma_1(x)\dfrac{\partial}{\partial x}u(\lambda,x) + (\sigma_2(x) \lambda + \gamma(x))u(\lambda,x) = 0
\]
to solutions $y(\lambda,x)=S(\lambda,x)u(\lambda,x)$ of the \textit{output}	LDE \eqref{eq:OutCC}
with the same spectral parameter:
\[
-\sigma_1(x)\dfrac{\partial}{\partial x} y(\lambda,x) + (\sigma_2(x) \lambda + \gamma_*(x))y(\lambda,x) = 0
\]
The first important result is Theorem \ref{thm:RealGen}, which relies on a realization
theorem of symmetric functions on Krein spaces \cite{bib:KreinReal}:
\begin{thm} \nonumber
Given a transfer function $S(\lambda,x)\in\Ir(\sigma_1(x), \sigma_2(x), \gamma(x), \gamma_*(x), \mathrm I)$
there exists a vessel
\[ \mathfrak{K_V} = (A, B(x), \mathbb X(x); \sigma_1(x), 
\sigma_2(x), \gamma(x), \gamma_*(x);
\mathcal{H}, \mathcal{E};\mathrm I_0),
\]
such that the transfer function of $\mathfrak{K_V}$ coincides with $S(\lambda,x)$, defined probably on a
smaller interval $\mathrm I_0\subseteq\mathrm I$.
\end{thm}
\noindent Since the existence of the inverse of $\mathbb X(x)$ plays so important role, we define in \eqref{eq:Tau}
$\tau = \det (\mathbb X^{-1}(0) \mathbb X(x))$. Since $\mathbb X(x)$ is a solution
of the Lyapunov equation \eqref{eq:XLyapunov}, this is a first sign, why we call this function
as "tau" function.

In order to even more emphasize the name and the role of the tau function, 
one have to consider SL equation \eqref{eq:SL}. 
In this case, one can uniquely reconstruct the potential from 
$S(\lambda,x_0)$ using the solution of Gelfand-Levitan equation 
\eqref{eq:GelfandLevitan}, constructed from $S(\lambda,x_0)$. More explicitly,
defining (see formulas \eqref{eq:OmeagaDef}, \eqref{eq:KDef})
\[ \begin{array}{ll}
\Omega(x,y) = \bbmatrix{1&0} B^*(x) \mathbb X^{-1}(x_0) B(y) \bbmatrix{1\\0}, \\
K(x,y) = -\bbmatrix{1&0} B^*(x) \mathbb X^{-1}(x) B(y) \bbmatrix{1\\0},
\end{array} \]
one finds that 
\[ K(x,y) + \Omega(x,y) + \int\limits_{x_0}^x K(x,t) \Omega(t,y) dt = 0,
\]
and the potential have the classical formula $q(x) = \dfrac{d}{dx} K(x,x) =
- 2 \dfrac{d^2}{dx^2} \ln \tau(x)$, which again explains the name for it. 

Basic ideas in this article come primarily from the work of M. Liv\v sic 
\cite{bib:Vortices}, which actually started in \cite{bib:defVess}. A generalization of Liv\v sic' vessel was developed in
\cite{bib:MyThesis, bib:MelVin1}, creating a comfortable background to learn linear differential equations with
a spectral parameter. It is important to notice that Liv\v sic' definition corresponds to dissipative vessels
(see Definition \ref{def:Vessel}) in our framework. In \cite{bib:SLVessels}, there is
presented an interesting research on finite dimensional vessels of the
equation \eqref{eq:SL}, which correspond to potentials, having purely discrete finite spectrum,
along with some interesting results related to differential algebras.

In section \ref{sec:SLDissImagine} there are discussed asymptotic behavior of vessel objects, in the case
the spectrum of $A$ is in $i\mathbb R^+$, which correspond to the classical case of being on the negative real line.
It turns out (Theorem \ref{thm:QBoundDiss}) that the potential of such a vessel satisfies 
$|q(x)| \leq \dfrac{Q}{(x-x_0)}$ for big enough $x$.

This theory would be less applicable without a concrete and simple example, which would show
how to construct a nontrivial example of a vessel. 
This task is successfully accomplished in Section
\ref{sec:ArbCurve}, where it is constructed a transfer function, having singularities (jumps)
along a given symmetric with respect to the imaginary axis curve (which may be unbounded).

Finally, in Sections \ref{sec:NLS}, \ref{sec:CanSys} we show that Non Linear Schr\" odinger equations and Canonical Systems fit in our framework.

\section{Vessels}
\subsection{\label{sec:DefVel}Definition of a vessel}
Before we define a notion of a vessel, one needs to define a list parameters, which will be fixed in many cases,
and thus is dealt separately. Then one defines a notion of a vessel, corresponding to these vessel parameters.
\begin{defn}
\label{def:VesPar}
Let $\sigma_1$, $\sigma_2$, $\gamma$, and $\gamma_*$  be
operators from a finite dimensional Hilbert space ${\mathcal E}$ to itself, 
locally integrable on an interval ${\mathrm I}=[a,b]$.
Suppose that $\sigma_1$ is differentiable and invertible on
${\mathrm I}$, and that the following relations hold:
\[ \begin{array}{lll}
\sigma_1(x) = \sigma_1^*(x), \quad \sigma_2(x) = \sigma_2^*(x) \\
\gamma(x) + \gamma(x)^* = \gamma_*(x) +
\gamma_*(x)^* = -\dfrac{d}{dx} \sigma_1(x),\quad x\in \mathrm I.
\end{array} \]
Then
$\sigma_1,\sigma_2, \gamma, \gamma_*$
and the interval $\mathrm I$ are called \textbf{vessel parameters} on $\mathcal E$.
\end{defn} \label{def:Vessel}
Before we define a notion of a vessel which involves an auxiliary Hilbert space $\mathcal H$ and operators
(for $x\in\mathrm I$)
\begin{equation} \label{eq:DefOperats} \begin{array}{llllllll}
A, \mathbb X(x) &:& \mathcal H &\rightarrow & \mathcal H, \\
B(x)						&:& \mathcal E &\rightarrow & \mathcal H
\end{array} \end{equation}
we have to consider some regularity assumptions.
We will assume that the operator $A$ may be unbounded with a domain $D(A)$.
Moreover, certain algebraic and differential relations 
will connect these operator, and as a result, we have to determine assumptions, which will ensure
that the relations between $A, \mathbb X(x), B(x)$ will become solvable equations.
\begin{defn}[Regularity assumptions] \label{def:UnOperts}
Operators $A, \mathbb X(x), B(x)$ are said to satisfy regularity assumptions on $\mathrm I$, if 
there exists a point $x_0\in\mathrm I$ such that 
\begin{enumerate}
	\item $B(x_0) \mathcal E \in D(A^n)$ for all $n\in\mathbb N$ and there exists $C>0$ such that
	\begin{equation} \label{eq:BIsSchwz}
		\| A^n B(x_0) \| \leq (C \sqrt{n})^n,
	\end{equation}
	\item The operator $\mathbb X(x)$ is self-adjoint and invertible for all $x\in\mathrm I$.
\end{enumerate}
\end{defn}
In order to show that such a requirement is fulfilled for some operators, we notice that
operator $A$ is usually isomorphic to the operator of multiplication by $t$ on $\mathbb R$.
Taking the initial condition $B(x_0) = e^{-t^2}$, we will obtain that
\[ \| A^{2n} B(x_0) \| = \int_{\mathbb R} t^{2n} e^{-t^2} dt = 
\dfrac{2}{2n+1} \| A^{2(n+1)} B(x_0) \|
\]
and the estimate above follows by induction.
This means that in the case $A$ is a multiplication by
$\mu$ on an unbounded curve $\Gamma$ and $\mathcal H = L^2(\Gamma)$, one can take 
$B(x_0)$ such that it decreases at infinity as $e^{-|\mu|^2}$.
For the vessel parameters one defines a notion of a vessel:
\begin{defn} \label{def:KV}
A \textbf{vessel} is a collection of operators and spaces
\begin{equation} \label{eq:DefKV}
\mathfrak{K_V} = (A, B(x), \mathbb X(x); \sigma_1(x), 
\sigma_2(x), \gamma(x), \gamma_*(x);
\mathcal{H}, \mathcal{E};\mathrm I), 
\end{equation}
where $\sigma_1(x), \sigma_2(x), \gamma(x), \gamma_*(x)$ and $\mathrm I$ are vessel parameters on $\mathcal E$. 
The spaces $\mathcal H$ is Hilbert and the operators $A, \mathbb X(x), B(x)$ are defined in \eqref{eq:DefOperats} so
that the regularity assumptions hold. Operators are subject to the following \textbf{vessel conditions}:
\begin{align}
\label{eq:DB} 0  =  \frac{d}{dx} (B(x)\sigma_1(x)) + A B(x) \sigma_2(x) + B(x) \gamma(x), \\
\label{eq:XLyapunov} A \mathbb X(x) + \mathbb X(x) A^*  +  B(x) \sigma_1(x) B^*(x) = 0, \\
\label{eq:DX} \frac{d}{dx} \mathbb X(x)  =  B(x) \sigma_2(x) B^*(x), \\
\label{eq:Linkage}
\gamma_*(x)  =  \gamma(x) + \sigma_2(x) B^*(x) \mathbb X^{-1}(x) B(x) \sigma_1(x) 
 - \sigma_1(x) B^*(x) \mathbb X^{-1}(x) B(x) \sigma_2(x).
\end{align}
\end{defn}
In order to understand why it is a well defined object, it is enough to show that
the equations, defining the vessel are solvable. We will see later that the equation
\eqref{eq:DB} is the key point of the construction and the rest will be easily constructed from
it.
\begin{thm}\label{thm:BConst}
Suppose that $B(x_0)$, A satisfy the condition \eqref{eq:BIsSchwz}, then
there exists a solution $B(x)$ of \eqref{eq:DB} with the value $B(x_0)$ for $x=x_0$.
Moreover, the estimate similar to \eqref{eq:BIsSchwz} holds
\[ \| A^n B(x) \| \leq (C(x) \sqrt{n})^n 
\]
\end{thm}
\noindent\textbf{Proof}
Before solve the equation \eqref{eq:DB}, notice that it is equivalent to
\[ B'(x) + A B(x) \sigma_2 \sigma_1^{-1} + B(x) [\gamma(x) + \sigma_1'(x)] \sigma_1^{-1} = 0.
\]
So, defining $E(x)$ such that
\[ E'(x) = [\gamma(x) + \sigma_1'(x)] \sigma_1^{-1} E(x), \quad E(x_0) = I,
\]
we will obtain the following equation
\begin{equation} \label{eq:BE}
 \dfrac{d}{dx}[B(x) E(x)] + A [B(x) E(x)] E^{-1}(x) \sigma_2 \sigma_1^{-1} E(x) = 0.
\end{equation}
Denote by $\Psi(x,\lambda)$ the solution of (substituting $B(x) E(x)$ with $\Psi(x,\lambda)$,
$A$ with $\lambda$, and denoting $E^{-1}(x) \sigma_2 \sigma_1^{-1} E(x)$ by $\widetilde E(x)$
in the last equation)
\[ \dfrac{d}{dx} \Psi(x,\lambda) + \lambda \Psi(x,\lambda) \widetilde E(x) = 0, \quad \Psi(x_0,\lambda) = I.
\]
From the Peano-Baker formula it follows that
\[ \begin{array}{lll}
\Psi(x,\lambda) & = I - \lambda\int\limits_{x_0}^{x} \widetilde E(y_1) dy_1 +
\lambda^2 \int\limits_{x_0}^x \int\limits_{x_0}^{y_1} \widetilde E(y_2) dy_2 
\widetilde E(y_1) dy_1  + \cdots = \\
& = \sum\limits_{n=0}^\infty \Psi_n(x) \lambda^n, \quad \Phi_0(x) = I,
\end{array} \]
and it is a well know result that the coefficient of this matrix satisfy the
relation $\Psi'_{n+1} = - \Psi_n(x) \widetilde E(x)$ and decrease as coefficients
of an exponential function $\|\Psi_n(x) \| \leq \dfrac{M}{n!}$. Let
\[ B_1(x) = \sum\limits_{n=0}^\infty A^n B(x_0) \Psi_n(x).
\]
Since $\|\Psi_n(x) \| \leq \dfrac{M}{n!}$ and by \eqref{eq:BIsSchwz} $\| A^n B(x_0) \| \leq (C \sqrt n)^n 
\leq C_1^n n!$, for some $C_1 < 1$,
we obtain that
\[ \| A^n B(x_0) \Psi_n(x) \| \leq M C_1^n
\]
which means that the series is absolutely convergent.
The same holds for the derivative. Differentiating this expression, we find that
\[ \begin{array}{lll}
\dfrac{d}{dx} B_1(x) & = 
\sum\limits_{n=0}^\infty A^n B(x_0) \dfrac{d}{dx} \Psi_n(x) = (\text{since }\Psi_0=I) \\
& =
-\sum\limits_{n=1}^\infty A^n B(x_0) \Psi_{n-1}(x) \widetilde E = 
-\sum\limits_{n=0}^\infty A^{n+1} B(x_0) \Psi_{n}(x) \widetilde E\\
& = - A \sum\limits_{n=0}^\infty A^n B(x_0) \Psi_{n}(x) \widetilde E = -A B_1(x) \widetilde E(x),
\end{array} \]
which means that $B_1(x)$ (and $B(x)E(x)$) satisfy the equation \eqref{eq:BE} with the
initial condition $B(x_0)$. Finally, $B(x) = B_1(x) E^{-1}(x)$ and this prove the existence.

Let us prove now the norm estimate for this function
\[ \begin{array}{lll}
\|A^n B(x)\| = \| \sum\limits_{k=0}^\infty A^{n+k} B(x_0) \Psi_k(x) \| \leq
\sum\limits_{k=0}^\infty (C\sqrt{n+k})^{n+k} \dfrac{M}{k!} =
C^n \sum\limits_{k=0}^\infty (\sqrt{n+k})^{n+k} \dfrac{M C^k}{k!} \\
\leq (C \sqrt{n})^n \sum\limits_{k=0}^\infty (\sqrt{1+\dfrac{k}{n}})^n (\sqrt{n+k})^k  \dfrac{M C^k}{k!}
\end{array} \]
Since $k!$ behaves like $\sqrt{2\pi k} (\dfrac{k}{e})^k$ asymptotically, there exists $C_1>0$ such that $\dfrac{1}{k!} \leq \dfrac{C_1 e^k}{k^k} $. Thus the last inequality becomes
\[ \begin{array}{lllll}
\|A^n B(x)\| & \leq (C \sqrt{n})^n \sum\limits_{k=0}^\infty (\sqrt{1+\dfrac{k}{n}})^n (\sqrt{n+k})^k  \dfrac{MC_1 (eC)^k}{k^k} \\
& \leq (C \sqrt{n})^n \sum\limits_{k=0}^\infty \dfrac{(\sqrt{1+\dfrac{k}{n}})^n}{k^{k/2}} (\sqrt{1+\dfrac{n}{k}})^k  MC_1 (eC)^k \\
& \leq (C \sqrt{n})^n \sum\limits_{k=0}^\infty \dfrac{(\sqrt{1+\dfrac{k}{n}})^n}{k^{k/2}} C_2 e^{n/2}  MC_1 (eC)^k \\
& \leq (C \sqrt{n})^n \sum\limits_{k=0}^\infty \dfrac{C_4 e^{k/2}}{k^{k/2}} C_2 e^{n/2}  MC_1 (eC)^k \\
& \leq (C \sqrt{n})^n C_5 C_4  C_2 e^{n/2}  MC_1 (eC)^k  \\
& \leq (C_6 \sqrt{n})^n, 
\end{array} \]
since the sum $\sum\limits_{k=0}^\infty \dfrac{e^{k/2}}{k^{k/2}}$ is finite.
\qed

The equation \eqref{eq:Linkage} is also called by M. Liv\v sic 
\cite{bib:defVess} as \textbf{linkage} condition. It turns out that the so called \textbf{Lyapunov equation} \eqref{eq:XLyapunov} is partially redundant.
\begin{lemma}[\textbf{Lyapunov condition permanence}] \label{lemma:Redund}
Suppose that $B(x)$ satisfies \eqref{eq:DB} and $\mathbb X(x)$ satisfies \eqref{eq:DX}, then if the Lyapunov equation
\eqref{eq:XLyapunov}
\[ A \mathbb X(x) + \mathbb X(x) A^*  + B(x) \sigma_1 B^*(x) = 0\]
holds for a fixed $x_0$, then it holds for all $x$. If $\mathbb X(x_0) = \mathbb X^*(x_0)$ then
$\mathbb X(x)$ is self-adjoint for all $x$.
\end{lemma}
\noindent\textbf{Proof:} By differentiating the Lyapunov equation, we will obtain
that LHS is constant. Since the derivative \eqref{eq:DX} $\dfrac{d}{dx}\mathbb X(x) = B(x) \sigma_2(x) B^*(x)$
is self-adjoint, $\mathbb X(x)$ will be self-adjoint, once $\mathbb X(x_0)$ is. \qed

We notice that $(\mathcal H, \mathbb X(x))$ form a Krein space, which is the same
as a set, but whose (Krein) inner product depends on $x$ and is differentiable. 
From the system theory \cite{bib:Brodskii, bib:Kalman} and the operator theory related to $J$-contractive functions
\cite{bib:bgr, bib:Potapov} we borrow some of the following additional characterizations of the vessel
\begin{defn} The vessel $\mathfrak{K_V}$ \eqref{eq:DefKV} is called
\begin{itemize}
	\item \textbf{dissipative}, if it is the case that $\mathbb X(x)>0$ for all values of $x\in\mathrm I$,
	\item \textbf{Pontryagin}, if $\mathbb X(x)$ has $\kappa\in\mathbb N$ negative squares at the
		right half plane for all values of $x\in\mathrm I$,
	\item \textbf{regular},  if all the operators $A, B(x), \mathbb X(x)$ are bounded operators for all $x$,
	\item \textbf{minimal}, if for all $x$ it holds that
		\begin{equation} \label{eq:minCond}
		 \operatorname{cl} \{ A^n B(x) \mathcal E\mid n\in\mathbb N \} = \mathcal H,
		\end{equation}
\end{itemize}
\end{defn}
\noindent where "$\operatorname{cl}$" stands for the closed span of the corresponding vectors.
One of the most important functions associated to the vessel is as follows \cite{bib:Brodskii}:
\begin{defn} The $\mathcal E\times\mathcal E$ valued function $S(\lambda, x)$ defined by
\begin{equation} \label{eq:S}
S(\lambda, x) = I - B^*(x) \mathbb X^{-1}(x) (\lambda I - A)^{-1} B(x) \sigma_1(x).
\end{equation}
is called the \textbf{transfer function} of the vessel $\mathfrak{K_V}$.
\end{defn}
It is extremely important and interesting case that the transfer functions and regular vessels are determined one
from the other as it will be shown in the following Section \ref{sec:RegVes} (it was first shown for regualar,
dissipative vessels in \cite{bib:MyThesis, bib:MelVinC}). 

\subsection{\label{sec:StandConstr}Standard construction of a vessel}
Let us show that one can easily construct vessels. For this to happen, choose two Hilbert spaces
$\mathcal H, \mathcal E$ and define three
operators $\mathbb X_0, A:\mathcal H \rightarrow \mathcal H$ and $B_0:\mathcal E\rightarrow\mathcal H$
such that $\mathbb X_0$ is invertible and the following equalities hold
\[ \mathbb X_0^* = \mathbb X_0, \quad
A \mathbb X_0 + \mathbb X_0 A^* +  B_0 \sigma_1(x_0) B_0^* = 0.
\]
Then solve \eqref{eq:DB} using Theorem \ref{thm:BConst}
\[ 0  =  \frac{d}{dx} (B(x)\sigma_1(x)) + A B(x)  \sigma_2(x) + B(x) \gamma(x), 
\quad B(x_0) =  B_0
\]
and solve the equation \eqref{eq:DX} by
\[ \mathbb X(x)  =  \mathbb X_0 + \int_{x_0}^ x B(y) \sigma_2(y) B(y)^* dy.
\]
Finally, define $\gamma_*(x)$ from $\gamma(x)$ using \eqref{eq:Linkage}. Thus a vessel is created (the interval 
$\mathrm I$ is defined in the proof):
\begin{lemma} The collection
\[ \mathfrak{K_V} = (A, B(x), \mathbb X(x); \sigma_1(x), 
\sigma_2(x), \gamma(x), \gamma_*(x);
\mathcal{H}, \mathcal{E};\mathrm I)
\]
is a vessel.
\end{lemma}
\noindent\textbf{Proof:} The equations \eqref{eq:DB}, \eqref{eq:DX}, and \eqref{eq:Linkage} are satisfied by the
construction. The Lyapunov equation \eqref{eq:XLyapunov} and the self-adjointness of $\mathbb X(x)$ follow from
Lemma \ref{lemma:Redund}. Since $\mathbb X_0$ is an invertible operator, there exists a non trivial interval
$\mathrm I$ (of length at least $\dfrac{1}{\|\mathbb X_0^{-1}\|}$) on which $\mathbb X^{-1}(x)$ exists.
\qed

We can obtain in this manner a rich family of vessels, since there exist standard models, creating operators
$\mathbb X_0, A, B_0$:
\begin{enumerate}
	\item Liv\v sic model of a non selfadjoint operator \cite{bib:BL}, where $\mathbb X=I$, $A+A^*+B J B^* = 0$, and
		$J$ is a signature matrix ($J=J^*=J^{-1}$),
	\item Theory of nodes, developed in \cite{bib:Brodskii},
	\item Krein space realizations for symmetric functions \cite{bib:KreinReal}, see the following Section
		\ref{sec:RegVes}.
	\item Vessels on curves (see Section \ref{sec:ArbCurve} in this article),
\end{enumerate}

\subsection{\label{sec:RegVes}Regular vessels}
We start from a realization theorem, which will enable us to construct regular vessels.
\subsubsection{A realization theorem for symmetric functions using Krein spaces}
Let $(\mathcal H,\langle\cdot,\cdot\rangle)$ be a Hilbert space. 
Let $\mathbb X$ be a self-adjoint \textit{invertible} operator on $\mathcal H$.
We define a sesquilinear form $[\cdot,\cdot]$ on $\mathcal H$ as $[u,v]=\langle \mathbb X u,v\rangle$. Define 
$\mathcal K$ to be as a set the same Hilbert space $\mathcal H$, but equipped with (indefinite) inner
product: $(\mathcal K,[\cdot,\cdot])$, which is called \textit{Krein space}. In the most general case one do not need
the invertability of the operator $\mathbb X$, but we will assume it for our purposes. For any operator $T$
on $\mathcal K$ we denote by $T^+$ he unique operator satisfying $[Tu,v]=[u,T^+v]$ for all $u,v\in\mathcal K$. Actually,
it follows that 
\[ T^+ = \mathbb X^{-1} T^* \mathbb X,
\]
where $T^*$ is the adjoint with respect to $\langle\cdot,\cdot\rangle$. 
The space $\mathcal H$ admits the decomposition
\[ \mathcal H = \mathcal H^+ \oplus \mathcal H^-
\]
such that $[u,u]>0$ for all $x\in\mathcal H^+$ and $[u,u]<0$ for all $x\in\mathcal H^-$. Moreover, the spaces \linebreak
$(\mathcal H^+,[\cdot,\cdot]), (\mathcal H^-,-[\cdot,\cdot])$ are complete with respect to the norms 
$[\cdot,\cdot]$ and $-[\cdot,\cdot]$ respectively.

Let $j=1,2$ and $(\mathcal K_j,[\cdot,\cdot]_j$ be Krein spaces and let $U_j$ be
linear operators in $\mathcal K_j$. The operator $U_1, U_2$ will be called \textit{weakly isomorphic},
if there exists dense subsets $L_j\subseteq\mathcal K_j$ such that $\Delta_j=L_j\cap D(U_j)$ (D - domain of) is dense
in $\mathcal K_j$ and $U_j(\Delta_j)\subseteq L_j$ and a bijection $V$ from $L_1$ to $L_2$ whic preserves
the indefinite inner product:
\[ [Vu,Vv]_1 = [v,u]_2, \quad \forall u,v \in L_1
\]
and has the properties $V\Delta_1=\Delta_2$, $VU_1u=U_2V u$ ($u\in\Delta_1$).

The following theorem
is taken from \cite[Theorem 3]{bib:KreinReal}, when we use a less powerful version of it
\begin{thm} \label{thm:RealizeSinK}
Let $S$ be a function, which is holomorphic on 
$\{ |\lambda|\mid\geq r\cup\{\infty\}$ with values in a Hilbert space $\mathcal E$.
Moreover, suppose that $S$ is symmetric with respect to the real axis:
\[ S(\lambda)=S^*(\bar \lambda)
\]
Then there exist a Krein space $\mathcal K$, a bounded self-adjoint operator $\widetilde A$ in $\mathcal K$
and $\widetilde \Gamma:\mathcal E\rightarrow \mathcal K$, such that
\[ S(\lambda) = S(\infty) + \widetilde\Gamma^+(\widetilde A-I\lambda)^{-1}\widetilde\Gamma.
\]
The space $\mathcal K$ can be chosen minimal:
\[ \mathcal K = \operatorname{cls}\{ (I\lambda - \widetilde A)^{-1} \widetilde\Gamma \mathcal E \mid \lambda\in \mathbb C\};
\]
then $\widetilde A$ is uniquely determined up to a weak isomorphism.
\end{thm}
Important remark, relevant to this research is that when the matrix $\mathbb X$ is strictly positive (and invertible)
we obtain the usual notion of the Hilbert space, equipped with the norm $[u,v] = \langle\mathbb Xu,v\rangle$.
\subsubsection{Regular vessels versus transfer functions}
Regular vessels, defined by bounded operators have a very good realization theory for their transfer functions.
Notice that in this case the functions $S(\lambda,x)$ are analytic at infinity (actually out of the spectrum of $A$)
with value $I$ there. It turns out that given just a transfer function itself one can reconstruct a vessel
using a theory of Krein realizations for functions, analytic at infinity (see Theorem \ref{thm:RealizeSinK}). 
Notice also that poles and singularities of $S(\lambda,x)$ with respect to $\lambda$ are determined 
by $A$ only and are independent of $x$.
Moreover, if a vessel is minimal (i.e. \eqref{eq:minCond} holds), standard theorems 
\cite{bib:bgr} in realization theory ensure that the singularities of $S(\lambda,x)$
occurs precisely at the spectrum of $A$.

In the next proposition, we summarize the properties of the transfer function of a regular vessel:
\begin{prop} 
\label{prop:PropS}
Let $\mathfrak{K_V}$ be a regular vessel and let $S(\lambda,x)$ be its transfer function. 
Then
\begin{enumerate}
	\item For all $x$, $S(\lambda, x)$ is an analytic function of $\lambda$
	in the neighborhood of $\infty$, where it satisfies: $S(\infty, x) = I$.
	\item For all $\lambda\not\in\SPEC(A)$, $S(\lambda, x)$ is a differentiable function of $x$.
	\item $S(\lambda,x)$ satisfies the symmetry condition 
	\begin{equation} \label{eq:Symmetry}
		S^*(-\bar\lambda,x) \sigma_1(x) S(\lambda,x) =  \sigma_1(x)
	\end{equation}
	for $\lambda$ in the domain of analyticity of $S(\lambda,x)$.
	\item Multiplication by $S(\lambda,x)$ maps solutions $u(\lambda,x)$ of the \textbf{input} 
	LDE with the spectral parameter $\lambda$:
	\begin{equation} \label{eq:InCC}
		-\sigma_1(x)\dfrac{\partial}{\partial x}u(\lambda,x) + (\sigma_2(x) \lambda + \gamma(x))u(\lambda,x) = 0
	\end{equation}
	to solutions $y(\lambda,x)=S(\lambda,x)u(\lambda,x)$ of the \textbf{output}	LDE with the same spectral parameter:
	\begin{equation} \label{eq:OutCC}
		-\sigma_1(x)\dfrac{\partial}{\partial x} y(\lambda,x) + (\sigma_2(x) \lambda + \gamma_*(x))y(\lambda,x) = 0
	\end{equation}
\end{enumerate}
\end{prop}
\noindent\textbf{Proof:}
These properties are easily checked, and follow from the definition
of $S(\lambda,x)$:
\[ S(\lambda, x) = I - B^*(x) \mathbb X^{-1}(x) (\lambda I - A)^{-1} B(x) \sigma_1(x). \]
The function $S(\lambda, x)$ is analytic
for $\lambda > \| A(x)\|$ and since all the operators are
bounded, we have $ S(\infty, x)
= I$. The second property follows from the differentiability assumptions
on the operators $\mathbb X(x), B(x)$. The third property follows from straightforward
calculations using the Lyapunov equation \eqref{eq:XLyapunov}:
\[ \begin{array}{lllll}
S(\mu, x)^* \sigma_1(x) S(\lambda, x) - \sigma_1(x) = \\
~~~~~~ - (\bar\mu+\lambda) \sigma_1(x)
B^*(x) (\bar\mu I - A^*)^{-1}\mathbb X^{-1}(x) (\lambda I - A)^{-1}  B(x) \sigma_1(x) = 0
\end{array} \]
for $\mu=-\bar\lambda$. The fourth property follows
directly from the definitions by plugging $y(\lambda,x) = S(\lambda, x) u(\lambda,x)$
into \eqref{eq:OutCC} and using \eqref{eq:InCC} for $u(\lambda,x)$, and the formula \eqref{eq:S} for $S(\lambda, x)$,
for which in turn we use vessel conditions in order to differentiate it. \qed

In fact, the converse of Proposition \ref{prop:PropS} holds. It was first proved in \cite{bib:MyThesis},
\cite[chapter 5]{bib:MelVinC} for the dissipative case (when $\mathbb X(x)>0$) and we
shall see later in Theorem \ref{thm:RealGen} that it holds for a regular vessel too.
We define the class of transfer functions, corresponding to the regular vessels as follows:
\begin{defn}[\cite{bib:MelVinC}] \label{def:IK}
The class
$\Ir=\Ir(\sigma_1(x), \sigma_2(x), \gamma(x), \gamma_*(x), \mathrm I)$
is the class consisting of functions $S(\lambda,x)$ of two variables possessing properties appearing in
Proposition \ref{prop:PropS}.
\end{defn}
Recall (see \cite{bib:CoddLev}) that to every LDE one can associate an
invertible matrix (or operator) function $\Phi(x,x_0)$, called
the {\sl fundamental solution}, which
takes value $I$ at some preassigned point $x_0$ and such that any other solution
$u(x)$ of the LDE, with initial condition $u(x_0) = u_0$ is of the form
\[ u(x) =
\Phi(x,x_0) u_0.
\]
Let  $\Phi(\lambda,x, x_0)$ and $\Phi_*(\lambda,x,
x_0)$ be the fundamental solutions of the input LDE \eqref{eq:InCC} and the
output LDE \eqref{eq:OutCC} respectively, where we have added in the
notation the dependence in $\lambda$. Then,
\begin{equation} \label{eq:SInttw}
S(\lambda, x) \Phi(\lambda,x,x_0) =
\Phi_*(\lambda,x,x_0) S(\lambda, x_0)
\end{equation}
and consequently $S(\lambda,x)$ satisfies the following LDE
\begin{equation}
\label{eq:DS}
\begin{array}{lll}
    \frac{\partial}{\partial x} S(\lambda,x)
    = \sigma_{1}^{-1}(x) (\sigma_{2}(x)
\lambda + \gamma_*(x)) S(\lambda,x)
-   S(\lambda,x)\sigma_1^{-1}(x) (\sigma_2(x) \lambda + \gamma(x)).
\end{array}
\end{equation}
The following properties of the fundamental matrices will be used in the sequel
\begin{prop} The following formulas hold
\begin{eqnarray} \label{eq:Phi*sig1Phi}
\Phi^*(x,-\bar\lambda) \sigma_1 \Phi_*(x,\lambda)  = \sigma_1, \\
\label{eq:DPhi^*S2Phi} \frac{\partial}{\partial x} [\Phi^*(x,\mu) \sigma_1 \Phi(x,\lambda)]  =
(\bar \mu + \lambda) \Phi^*(x,\mu) \sigma_2 \Phi(x,\lambda).
\end{eqnarray}
and the same formulas hold, substituting $\Phi$ by $\Phi_*$.
\end{prop}
\noindent\textbf{Proof:} Immediate from the definitions. \qed

The next theorem shows that the class $\Ir$ is well-defined in the sense, that given a function, one can also find a 
corresponding to it vessel.
\begin{thm} \label{thm:RealGen}
Given a transfer function $S(\lambda,x)\in\Ir(\sigma_1(x), \sigma_2(x), \gamma(x), \gamma_*(x), \mathrm I)$,
there exists a regular vessel
\[ \mathfrak{K_V} = (A, B(x), \mathbb X(x); \sigma_1(x), 
\sigma_2(x), \gamma(x), \gamma_*(x);
\mathcal{H}, \mathcal{E};\mathrm I_0),
\]
such that the transfer function of $\mathfrak{K_V}$ coincides with $S(\lambda,x)$, defined probably on a
smaller interval $\mathrm I_0\subseteq\mathrm I$.
\end{thm}
\noindent\textbf{Proof:} Fix a point $x_0\in\mathrm I$ and define a function $Q(\lambda)$ using Calley transform,
which satisfies
\begin{equation} \label{eq:SfromQ}
S(-i\lambda,x_0) = ( I + \dfrac{i}{2} Q(\lambda) \sigma_1(x_0)) (I - \dfrac{i}{2} Q(\lambda)\sigma_1(x_0))^{-1}.
\end{equation}
Actually, this function is given by 
\[ Q(\lambda) = 2 i \sigma_1^{-1}(x_0) (I - S(-i\lambda,x_0)) (I + S(-i\lambda,x_0))^{-1}
\]
and is well-defined at the neighborhood of infinity with value $0$ there. 
Then from the equality
\[
(I-S^*(-i\bar\lambda,x_0))\sigma_1^{-1}(x_0)(I + S(-i\lambda,x_0)) =
 - (I+S^*(-i\bar\lambda,x_0))\sigma_1^{-1}(x_0)(I - S(-i\lambda,x_0)),
\]
resulting from the symmetry condition \eqref{eq:Symmetry}, considered with $-i\lambda$ instead of $\lambda$,
it follows that $Q(\lambda)^* = Q(\bar\lambda)$ and $Q(\lambda)$ is zero at the neighborhood of
$\lambda=\infty$. Thus from Theorem \ref{thm:RealizeSinK} it follows that $Q(\lambda,x)$ admits a
following Krein space realization
\[ Q(\lambda) = \Gamma^+ (\widetilde A - \lambda I)^{-1} \Gamma
\]
for a Krein space $\widetilde{\mathcal K}$ with self-adjoint operator $\widetilde A$ in $\widetilde{\mathcal K}$, 
and $\Gamma\in\mathcal E\rightarrow \mathcal K$. Inserting further this realization formula into
\eqref{eq:SfromQ} and simplifying we obtain:
\[ \begin{array}{lllllll}
S(-i\lambda,x_0) & = ( I + \dfrac{i}{2} Q(\lambda) \sigma_1(x_0)) (I - \dfrac{i}{2} Q(\lambda)\sigma_1(x_0))^{-1} = \\
& = ( 2I - I + \dfrac{i}{2} Q(\lambda) \sigma_1(x_0)) (I - \dfrac{i}{2} Q(\lambda)\sigma_1(x_0))^{-1}  =
 -I + 2 (I - \dfrac{i}{2} Q(\lambda)\sigma_1(x_0))^{-1}  = \\
& = -I + 2 (I - \dfrac{i}{2} \Gamma^+ (\widetilde A - \lambda I)^{-1} \Gamma \sigma_1(x_0))^{-1}
\end{array} \]
There is a simple formula \cite{bib:bgr} for evaluating the inverse of a matrix in a realized form:
\[ (I - \dfrac{i}{2} \Gamma^+ (\widetilde A - \lambda I)^{-1} \Gamma \sigma_1(x_0))^{-1} =
I + \dfrac{i}{2} \Gamma^+ (\widetilde A^\times - \lambda I)^{-1} \Gamma \sigma_1(x_0),
\]
where $\widetilde A^\times = \widetilde A - \dfrac{i}{2} \Gamma \sigma_1(x_0) \Gamma^+$. So, the last formula becomes
\begin{equation} \label{eq:SPreKrein} \begin{array}{lllllll}
S(-i\lambda,x) & = -I + 2 (I - i \Gamma^+ (\widetilde A - \lambda I)^{-1} \Gamma \sigma_1(x_0))^{-1} = \\
& = -I + 2(I+\dfrac{i}{2} \Gamma^+ (\widetilde A^\times - \lambda I)^{-1} \Gamma \sigma_1(x_0)) = 
 I + i \Gamma^+ (\widetilde A^\times - \lambda I)^{-1} \Gamma \sigma_1(x_0) = \\
& = I -  \Gamma^+ (i \widetilde A^\times - i\lambda I)^{-1} \Gamma \sigma_1(x_0) 
\end{array} \end{equation}
Let us define $A = - i \widetilde A^\times$ then we obtain that
\begin{equation} \label{eq:LyapPre}
 \begin{array}{lllllll}
A + A^{+} & = -i \widetilde A^\times + i (\widetilde A^\times)^+ = \\
& = -i(\widetilde A - \dfrac{i}{2} \Gamma \sigma_1(x_0) \Gamma^+) + i (\widetilde A^+ + \dfrac{i}{2} \Gamma \sigma_1(x_0) \Gamma^+) = \\
& =  -\Gamma \sigma_1(x_0) \Gamma^+,
\end{array} \end{equation}
since $\widetilde A$ is selfadjoint. Using the formulas for Krein space adjoint
\[ A^{+} = \widetilde{\mathbb X}^{-1} A^* \widetilde{\mathbb X}, \quad
\Gamma^+ = \widetilde{\mathbb X} \Gamma^*
\]
the last equation \eqref{eq:LyapPre} is
\[ \begin{array}{lllllll}
A +  A^{+} = - \Gamma \sigma_1(x_0) \Gamma^+ \Leftrightarrow \\
A + \widetilde{\mathbb X}^{-1} A^* \widetilde{\mathbb X} =
 \Gamma \sigma_1(x_0) \Gamma^* \widetilde{\mathbb X}  \Leftrightarrow \\
A \widetilde{\mathbb X}^{-1} + \widetilde{\mathbb X}^{-1} A^* = \Gamma \sigma_1(x_0) \Gamma^* 
\end{array}  \]
which is exactly the Lyapunov equation \eqref{eq:XLyapunov} at $x_0$ after defining
$B_0 = \Gamma$ and $\mathbb X_0 = \widetilde{\mathbb X}^{-1}$. Thus we obtain that
\begin{eqnarray} \label{eq:St20Krein}
S(\lambda,x_0) = I - B_0^* \mathbb X_0^{-1} (\lambda I - A)^{-1} B_0 \sigma_1(x_0), \\
A \mathbb X_0 + \mathbb X_0 A^* + B_0 \sigma_1(x_0) B^*_0 = 0, \quad \mathbb X_0 = \mathbb X_0^*.
\end{eqnarray}
As a result, we can use the standard construction of a vessel (see Section \ref{sec:StandConstr}), starting from $A, \mathbb X_0, B_0$
and to obtain a vessel, whose transfer function $Y(x,s)$ maps solution of the input LDE \eqref{eq:InCC} to solution
of the output LDE \eqref{eq:OutCC} for some $\gamma'_*(x)$.
Thus the two functions $S(\lambda,x)$ and $Y(\lambda,x)$
have the same value at $x_0$ and map solutions of the same input LDE to (possibly different) output LDEs:
\[ \begin{array}{ll}
S(\lambda,x) =
\Phi_*(\lambda,x,x_0) S(\lambda,x_0) \Phi^{-1}(\lambda,x,x_0), \\
Y(\lambda,x) = \Phi'_*(\lambda,x,x_0)
S(\lambda,x_0) \Phi^{-1}(\lambda,x,x_0).
\end{array} \]
Consequently, the function $S^{-1}(\lambda,x)Y(\lambda,x)$ is equal to
$I$ at infinity and is entire. By Liouville's theorem it is a
constant function and is equal to $I$. Thus
\[
\Phi_*(\lambda,x,x_0) = \Phi'_*(\lambda,x,x_0),
\]
from where we obtain that
\[
\Phi_*^{-1}(\lambda,x,x_0) \Phi'_*(\lambda,x,x_0) = I.
\]
Differentiating both sides of this last equation we are led to
\[ \begin{array}{ll}
0 & = \dfrac{\partial}{\partial x} [\Phi_*^{-1}(\lambda,x,x_0)
\Phi'_*(\lambda,x,x_0)] = \\
& = \Phi_*^{-1}(\lambda,x,x_0) \sigma_1^{-1}(x) (-\gamma(x) +
\gamma'(x)) \Phi'_*(\lambda,x,x_0).
\end{array}  \]
Since the matrices $\Phi_*(\lambda,x,x_0),
\Phi'_*(\lambda,x,x_0),\sigma_1(x) $ are invertible we obtain
that $\gamma_*(x) = \gamma_*'(x)$. It is remained to notice that constructed $\mathbb X(x), B(x), \gamma_*(x)$
satisfy the vessel conditions and as a result we obtain that the collection 
\[ \mathfrak{K_V} = (A, B(x), \mathbb X(x); \sigma_1(x), 
\sigma_2(x), \gamma(x), \gamma_*(x);
\mathcal{H}, \mathcal{E};\mathrm I_0),
\]
is a vessel whose characteristic function $Y(\lambda,x)$ coincides with $S(\lambda,x)$ on $\mathrm I_0$. \qed

Finally, for creating a complete picture of the correspondence between vessels and their transfer function, we
have to recall the following theorem
\begin{thm} Supose that we are given two regular \textbf{minimal} vessels
\[ \begin{array}{lll}
\mathfrak{K_V} = (A, B(x), \mathbb X(x); \sigma_1(x), 
\sigma_2(x), \gamma(x), \gamma_*(x);
\mathcal{H}, \mathcal{E};\mathrm I) \\
\widetilde{\mathfrak{K_V}} = (\widetilde A, \widetilde B(x), \widetilde{\mathbb X}(x); \sigma_1(x), 
\sigma_2(x), \gamma(x), \gamma_*(x);
\widetilde{\mathcal{H}}, \mathcal{E};\mathrm I)
\end{array} \]
then the transfer functions of these vessels are equal at the neighborhood of infinity if and only if
there exists an invertible densely definable operator
\[ T: \mathcal H \rightarrow \widetilde{\mathcal{H}},
\]
such that
\[ \widetilde A = T A T^{-1}, \quad 
\widetilde B(x) = T B(x), \quad
\widetilde{\mathbb X}(x) = V \mathbb X(x) V^*.
\]
\end{thm}
\noindent\textbf{Proof(outline):}
One of the directions is simple. Suppose that there exists such an operator $T$, then the transfer function of
$\widetilde{\mathfrak{K_V}}$ is
\[ \begin{array}{lllll}
\widetilde S(\lambda,x) & = 
I -  \widetilde B^*(x) \widetilde{\mathbb X}^{-1}(x) (\lambda I - \widetilde A)^{-1} B(x) \sigma_1(x) \\
& =  I - B^*(x) T^* (T \mathbb X(x) T^*)^{-1} (\lambda I - T A T^{-1})^{-1} T B(x) \sigma_1(x)  \\
& =  I - B^*(x) (\mathbb X(x))^{-1} (\lambda I -  A )^{-1} B(x) \sigma_1(x)  \\
& = S(\lambda,x).
\end{array} \]
For the converse direction, one uses the fact that this theorem holds when we fix $x_0\in\mathrm I$
(this a standard theorem in the realization theory of functions \cite{bib:helton, bib:bgr}). Then using the idea
appearing in Theorem \ref{thm:RealGen} (see last argument) that if two transfer functions are identical at $x_0$,
then they are identical for all $x$, we will obtain the desired result.
\qed 

\noindent\textbf{Remark:} The second part of this theorem, we could probably call as the \textbf{permanence of the
similarity operator}, in the spirit of Lemmas \ref{lemma:Redund}, \ref{lemma:detS} and 
Theorems \ref{thm:PermMin}, \ref{thm:S=S-1*}.

\subsection{Additional properties of the transfer function of a vessel}
The idea of this realization theorem may lead the reader to a conclusion that fixing 
$\sigma_1, \sigma_2, \gamma$ and varying the initial data $S(\lambda,x_0)$ one should obtain different
$\gamma_*$, uniquely determined by the vessel condition \eqref{eq:Linkage}. Unfortunately, it is not true
\begin{prop} \label{prop:SY}
Suppose that there exists a symmetric function $Y(\lambda)$, which commutes
with $\Phi(\lambda,x,x_0)$ and suppose that a function $S(\lambda,x)$ corresponds to
vessel parameters $\sigma_1,\sigma_2,\gamma,\gamma_*$. Then the function 
$S(\lambda,x) Y(\lambda)$ corresponds to the same vessel parameters $\sigma_1,\sigma_2,\gamma,\gamma_*$.
\end{prop}
\noindent\textbf{Proof:} Using formula \eqref{eq:SInttw} we obtain that
\[ S(\lambda, x) =
\Phi_*(\lambda,x,x_0) S(\lambda, x_0) \Phi^{-1}(\lambda,x,x_0) .
\]
Consequently,
\[ S(\lambda, x) Y(\lambda) = 
\Phi_*(\lambda,x,x_0) S(\lambda, x_0) Y(\lambda) \Phi^{-1}(\lambda,x,x_0) 
\]
intertwines solutions of the input \eqref{eq:InCC} and the output \eqref{eq:OutCC}
ODEs with the spectral parameter $\lambda$,
and is identity at infinity, because $S(\lambda, x)$ and $Y(\lambda)$ and their product are such. 
Symmetry is easy to check and by the definition the function $S(\lambda, x) Y(\lambda)$
corresponds to the same vessel parameters as $S(\lambda, x)$. \qed

Here are some interesting properties of the determinant of transfer functions, which will be used later
\begin{lemma}[\textbf{Permanence of $\det S$}] \label{lemma:detS} The determinant of the matrix-function
$S(\lambda,x)$ is $x$-independent and it holds that
\[ \det S(\lambda,x) = \det S(\lambda,x_0), \quad \lambda\not\in\SPEC{A}.
\]
For $\lambda$ on the imaginary axis, it holds that $|\det S(\lambda,x_0)| = 1, \lambda\not\in\SPEC{A}$.
\end{lemma}
\noindent\textbf{Proof:} Using equation \eqref{eq:DS} we obtain that
\[ \begin{array}{llll}
\dfrac{\dfrac{\partial}{\partial x} \det S(\lambda,x)}{\det S(\lambda,x)} =
\TR (S^{-1}(\lambda,x)\dfrac{\partial}{\partial x} S(\lambda,x)) = \\
= \TR(S^{-1}(\lambda,x)
\sigma_{1}^{-1}(x) (\sigma_2(x) \lambda + \gamma_*(x)) S(\lambda,x)
-  S(\lambda,x)\sigma_1^{-1}(x) (\sigma_2(x) \lambda + \gamma(x)) = \\
= \TR(\sigma_{1}^{-1}(x) (\sigma_2(x) \lambda + \gamma_*(x)) 
- \sigma_1^{-1}(x) (\sigma_2(x) \lambda + \gamma(x)) = \TR(\sigma_{1}^{-1}(x) (\gamma_*(x) - \gamma(x)) = \\
= \TR(\sigma_1^{-1}(x)B^*(x)\mathbb X^{-1}(x) B(x) \sigma_2(x) -
\sigma_1^{-1}(x)\sigma_2(x) B^*(x)\mathbb X^{-1}(x) B(x) \sigma_1(x) = 0,
\end{array} \]
we have used the property of trace $\TR(AB)=\TR(BA)$ and the linkage condition \eqref{eq:Linkage}. So, \linebreak
$\dfrac{\partial}{\partial x} \det S(\lambda,x) = 0$ and the result follows.

The second part of the lemma follows from taking determinant at the symmetry condition \eqref{eq:Symmetry}
and the fact that $-\bar\lambda=\lambda$, if $\lambda$ is on the imaginary axis.
\qed

Another interesting property of transfer functions is that when $\mathbb X(x) > 0$, they define a Kernel of
a Reproducing Kernel Hilbert Space \cite{bib:AronsjKErnels} as the following lemma states
\begin{lemma} \label{lemma:PosKernels}
Suppose that $\mathbb X(x) > 0$, then the Kernels
\[ K_1(\lambda,\mu,x) = \dfrac{\sigma_1(x) - S^*(\mu,x) \sigma_1(x) S(\lambda,x)}{\bar\mu + \lambda},
\quad K_2(\lambda,\mu,x) = \dfrac{\sigma_1^{-1}(x) - S(\lambda,x) \sigma_1^{-1}(x) S^*(\mu,x)}{\bar\mu + \lambda},
\]
are positive.
\end{lemma}
\noindent\textbf{Proof:} Using Lyapunov equation \eqref{eq:XLyapunov}, one obtains that
\begin{equation} \label{eq:K1}
K_1(\lambda,\mu,x) =
\sigma_1(x) B^*(x) (\bar\mu I - A^*)^{-1}\mathbb X^{-1}(x) (\lambda I - A)^{-1}  B(x) \sigma_1(x),
\end{equation}
and
\[
K_2(\lambda,\mu,x) = 
 B^*(x) \mathbb X^{-1}(x) (\bar\mu I - A^*)^{-1}\mathbb X(x) (\lambda I - A)^{-1} \mathbb X^{-1}(x) B(x),
\]
from where the positivity follows.
\qed

M. Livsi\v c model of a non self-adjoint operator \cite{bib:BL} uses (implicitly) the
assumption $\mathbb X_0=I$ and corresponds to a dissipative vessel
by the definition. The next theorem shows that an arbitrary dissipative vessel can be brought to a Livsi\v c
form for a fixed value $x_0$. Then one can use the standard construction, based on the Livsi\v c model at $x_0$.
\begin{thm} \label{thm:X0isI}
For a dissipative vessel there exists a Hilbert space similarity $V:\mathcal H\rightarrow \mathcal H$,
so that the new vessel
\[ \widetilde{\mathfrak{K_V}} = (\widetilde A, \widetilde B(x), \widetilde{\mathbb X}(x); \sigma_1(x), 
\sigma_2(x), \gamma(x), \gamma_*(x);
\mathcal{H}, \mathcal{E};\mathrm I_0)
\]
where 
\[ \widetilde A = V^{-1} A V, \quad \widetilde B(x) = V^{-1} B(x), 
\quad \widetilde{\mathbb X}(x) = V^{-1} \mathbb X(x) V^{-1},
\]
has the same transfer function as $\mathfrak{K_V}$ but satisfies additionally
\[ \widetilde{\mathbb X}(x_0) = I.
\]
\end{thm}
\noindent\textbf{Proof:} Define $V = \sqrt{\mathbb X(x_0)} > 0$, which exists, since $\mathbb X(x_0) > 0$. Then check that
the transfer function of $\widetilde{\mathfrak{K_V}}$ coincides with that of $\mathfrak{K_V}$:
\[  \begin{array}{llll}
\widetilde S(\lambda,x) & = 
I - \widetilde B^*(x) \widetilde{\mathbb X}^{-1}(x) (\lambda I - \widetilde A)^{-1} \widetilde B(x) \sigma_1 = \\
& = I - B^*(x) V^{-1} ( V^{-1} \mathbb X(x) V^{-1})^{-1} (\lambda I - V^{-1} A V)^{-1} V^{-1} B(x) \sigma_1 =
S(\lambda,x),
\end{array} \]
after cancellations. Finally notice that
\[ \widetilde{\mathbb X}(x_0) = V^{-1} \mathbb X(x_0) V^{-1} = I.
\]
\qed

The minimality condition turns out to be independent of $x$ in the sense that if it holds for one $x_0$,
then it holds for all $x$:
\begin{thm}[\textbf{Permanence of minimality}]\label{thm:PermMin} Suppose that we are given a vessel $\mathfrak{K_V}$, which is minimal at $x_0\in\mathrm I$.
Then the vessel is minimal for all $x\in\mathrm I$.
\end{thm}
\noindent\textbf{Proof:} Let us show it for the regular case first.
Suppose that the realization is minimal \eqref{eq:minCond}: 
$\operatorname{cl} \{ A^n B(x_0) \mathcal E\mid n\in\mathbb N \} = \mathcal H$. Using regularity assumption \eqref{eq:BIsSchwz},
we can represent $B(x)$ using a fundamental matrix $\Phi_1(\lambda,x,x_0)$:
\[ B(x) = \dfrac{1}{2\pi i} \oint (\lambda I - A)^{-1} B(x_0) \Phi_1(\lambda,x,x_0) d\lambda.
\]
Since $\Phi_1(\lambda,x,x_0)$ and its inverse are entire, we can use Taylor series in $\lambda$
\[ \Phi_1^{-1}(\lambda,x,x_0) = \sum \phi_k(x) \lambda^k.
\]
Finally, applying $A^{n+k} B(x)$ to matrices $\phi_k(x)$ and taking sums, whose convergence follows from the
analyticity of $\Phi_1^{-1}(\lambda,x,x_0)$ we will obtain that
\[  \begin{array}{llll}
\sum\limits_{k=0}^\infty A^{n+k} B(x)  \phi_k(x) = 
A^n \dfrac{1}{2\pi i} \sum\limits_{k=0}^\infty A^k (\lambda I - A)^{-1} B(x_0) \Phi(\lambda,x,x_0)\phi_k(x) = \\
= \dfrac{1}{2\pi i} \oint (\lambda I - A)^{-1} B(x_0) \Phi(\lambda,x,x_0) \sum_k \lambda^k \phi_k(x) d\lambda =
A^n \dfrac{1}{2\pi i} \oint (\lambda I - A)^{-1} B(x_0) d\lambda =
A^n B_0
\end{array} \]
and so 
\[
\operatorname{cl} \{ A^n B(x_0) \mathcal E\mid n\in\mathbb N \} \subseteq
\operatorname{cl} \{ A^n B(x) \mathcal E\mid n\in\mathbb N \},
\]
from where the minimality for all $x$ follows.

In the unbounded case, immediate consequence of the formula 
\[ B(x) = \sum A^k B_0 \phi_k(x) E^{-1}(x)
\] 
is that
\[ A^n B(x) \mathcal E \subseteq \operatorname{cl} \{ A^n B_0 \mathcal E\}.
\]
On the other hand, solving the differential equation \eqref{eq:DB} starting from
$x$ down to $x_0$, we will obtain that similar formulas hold with $x$ and $x_0$ interchanged
and in this case
\[ A^n B(x_0) \mathcal E \subseteq \operatorname{cl} \{ A^n B(x) \mathcal E\},
\]
from where it follows the permanence of the minimality. \qed

The next theorem shows that the symmetry condition (as the Lyapunov equation \eqref{eq:XLyapunov})
can be checked in one point.
This theorem is similar to the property of a solution of a Riccati
equation \cite[theorem 2.1]{bib:Zelikin}
\begin{thm}[\textbf{Permanence of symmetry}] \label{thm:S=S-1*} Suppose that the vessel parameters
satisfy Livsic (not M. Livsi\v c) condition.
Suppose that $S(\lambda,x)$ is a differentiable function of $x$ for each $\lambda$,
analytic in $\lambda$ for each $x$, 
except for a set of singular points, and satisfies $S(\infty,x) = I$. Suppose also that
$S(\lambda,x)$ is an intertwining function of LDEs \eqref{eq:InCC} and \eqref{eq:OutCC}.
Then if the symmetry condition \eqref{eq:Symmetry}
\[
S(\lambda,x) = \sigma_1^{-1}(x) S^{-1*}(-\bar\lambda,x) \sigma_1(x)
\]
holds for $x_0$, then it holds for all values of $x$.
\end{thm}
\noindent\textbf{Proof:} Since $S(\lambda,x)$ intertwines solutions of \eqref{eq:InCC} and \eqref{eq:OutCC}, then
it satisfies the differential equation \eqref{eq:DS}
\begin{eqnarray*} \frac{\partial}{\partial x} S(\lambda,x)
    = \sigma_{1}^{-1}(x) (\sigma_{2}(x)
\lambda + \gamma_*(x)) S(\lambda,x)-\\
-   S(\lambda,x)\sigma_1^{-1}(x) (\sigma_2(x) \lambda + \gamma(x)).
\end{eqnarray*}
Consequently, using properties of $\gamma_*, \gamma$ appearing in Definition \ref{def:VesPar}
we obtain that the function 
$\sigma_1^{-1}(x) S^{-1*}(-\bar\lambda,x) \sigma_1(x)$ satisfies the same differential equation.
If these two functions are equal at $x_0$, from the uniqueness of solution for a differential
equation with continuous coefficients, they are also equal for all $x$.
\qed

\subsection{The tau function of a vessel}
Following the ideas presented in \cite{bib:SLVessels} we define the tau function of the vessel 
\ref{def:KV} in the following way
\begin{defn} \label{def:Tau} For a given vessel (see Definition \ref{def:KV})
\[
\mathfrak{K_V} = (A, B(x), \mathbb X(x); \sigma_1(x), 
\sigma_2(x), \gamma(x), \gamma_*(x);
\mathcal{H}, \mathcal{E};\mathrm I)
\]
the tau function $\tau(x)$ is defined as
\begin{equation} \label{eq:Tau} \tau = \det (\mathbb X^{-1}(x_0) \mathbb X(x))
\end{equation}
for an arbitrary point $x_0\in\mathrm I$.
\end{defn}
Notice that using vessel condition \eqref{eq:DX} $\mathbb X(x)$ has the formula
\[ \mathbb X(x) = \mathbb X(x_0) + \int\limits_{x_0}^x B^*(y) \sigma_2 B(y) dy,
\]
and as a result
\[ \mathbb X^{-1}(x_0) \mathbb X(x) = I + \mathbb X^{-1}(x_0) \int\limits_{x_0}^x B^*(y) \sigma_2 B(y) dy.
\]
Since $\sigma_2$ has finite rank for $\dim \mathcal E<\infty$, 
this expression is of the form $I + T$, for a trace class operator $T$ and since 
$\mathbb X_0$ is an invertible operator, there exists a non trivial interval (of length at least $\dfrac{1}{\|\mathbb X_0^{-1}\|}$) on which $\mathbb X(x)$ and $\tau(x)$ are defined. Recall \cite{bib:GKintro} that a function $F(x)$ from (a, b) into the group G (the set of bounded invertible operators on H of the form I + T, for
a trace-class operator $T$) is said to be differentiable if $F(x) -I$ is \textit{differentiable} as a map into the trace-class operators. In our case,
\[ \dfrac{d}{dx} (\mathbb X^{-1}(x_0)\mathbb X(x)) = 
\mathbb X^{-1}(x_0) \dfrac{d}{dx} \mathbb X(x) =
\mathbb X^{-1}(x_0) B(x)\sigma_2B^*(x)
\]
exists in trace-class norm. Israel Gohberg and Mark Krein \cite[formula 1.14 on p. 163]{bib:GKintro}
proved that if $\mathbb X^{-1}(x_0)\mathbb X(x)$ is a differentiable function
into G, then $\tau(x) = \SP (\mathbb X^{-1}(x_0)\mathbb X(x))$
\footnote{$\SP$ - stands for the trace in the infinite dimensional space.} is a differentiable map into $\mathbb C^*$ with
\begin{multline} \label{eq:GKform}
\dfrac{\tau'}{\tau}  = \SP (\big(\mathbb X^{-1}(x_0)\mathbb X(x)\big)^{-1} 
\dfrac{d}{dx} \big(\mathbb X^{-1}(x_0)\mathbb X(x)\big)) = \SP (\mathbb X(x)' \mathbb X^{-1}(x)) = \\
= \SP (B(x)\sigma_2 B^*(x) \mathbb X^{-1}(x)) =
\TR (\sigma_2 B^*(x) \mathbb X^{-1}(x)B(x)).
\end{multline}
Since any two realizations of a symmetric function
are (weakly) isomorphic, one obtains from standard theorems \cite{bib:bgr} 
in realization theory of analytic at infinity functions that they will
have the same tau function up to a scalar, 
i.e. this notion is independent of the realization we choose for the given function $S(\lambda,x)$.
And we can call this property as the \textbf{permanence of the tau function}.

\section{Sturm-Liouville vessels}
Now we are ready to consider a particular example of vessels, which corresponds to the Sturm-Liouville differential
equation \eqref{eq:SL}. Some definitions in the general theory of vessels (such as the tau function) are
inspired by this particular example.
In order to obtain a SL vessel we choose 
$\mathcal E = \mathbb C^2$, i.e., a Hilbert space of dimension $2$ and make the following
\begin{defn} \label{def:SLparam}
The Sturm Liouville (SL) vessel parameters are defined to be
\[ \sigma_1 = \bbmatrix{0 & 1 \\ 1 & 0},
\sigma_2 = \bbmatrix{1 & 0 \\ 0 & 0},
\gamma = \bbmatrix{0 & 0 \\ 0 & i},
\gamma_*(x)=\bbmatrix{-i (\beta'(x) - \beta^2(x)) & -\beta(x) \\ \beta(x) & i}
\]
for a real valued function $\beta(x)$ differentiable an interval $\mathrm I$.
\end{defn}
Suppose now that we have a vessel $\mathfrak{K_V}$ realizing these vessel parameters. Then it turns out that
multiplication by the transfer function maps solutions of the trivial SL equation ($q(x)=0$) to a more complicated one.
As a result it can be considered as a B\" acklund transformation \cite{bib:CoddLev}. One can check that the Crum
transformations \cite{bib:Crum} are a particular case of vessels constructions 
(see \cite[Section 3.2]{bib:SLVessels} for details)

Denoting $u(\lambda,x) = \bbmatrix{u_1(\lambda,x)\\ u_2(\lambda,x)}$ we shall obtain that the input compatibility
condition \eqref{eq:InCC} is equivalent to
\[ \left\{ \begin{array}{lll}
-\frac{\partial^2}{\partial x^2} u_1(\lambda,x) = -i\lambda u_1(\lambda,x), \\
u_2(\lambda,x) = - i \frac{\partial}{\partial x} u_1(\lambda,x).
\end{array}\right.
\]
The output $y(\lambda,x) = \bbmatrix{y_1(\lambda,x)\\y_2(\lambda,x)} = S(\lambda,x) u(\lambda,x)$ satisfies the output
equation \eqref{eq:OutCC}, which is equivalent to
\[ \left\{ \begin{array}{lll}
-\frac{\partial^2}{\partial x^2} y_1(\lambda,x) + 2 \beta'(x) y_1(\lambda,x) = -i\lambda y_1(\lambda,x), \\
y_2(\lambda,x) = - i [ \frac{\partial}{\partial x} y_1(\lambda,x) - \beta(x) y_1(\lambda,x)].
\end{array}\right.
\]
Observing the first coordinates $u_1, y_1$ of the vector-functions $u,y$ we can see that multiplication by
$S(\lambda,x)$ maps solution of the trivial SL equation to solutions of the more complicated one,
defined by the potential 
\begin{equation} \label{eq:qbeta}
q(x) = 2 \beta'(x).
\end{equation}
If we denote the fundamental matrix of the input SL equation \eqref{eq:InCC} as following
\begin{equation} \label{eq:Phi}
 \Phi(x,\lambda) = \bbmatrix{\cos(s x) & \dfrac{i \sin(s x)}{s} \\  i s \sin(s x) & \cos(s x)}, \quad s^2 = -i\lambda
\end{equation}
then $S(x,\lambda) \Phi(x,s)$ is the fundamental matrix for solutions of the output SL equation
\eqref{eq:OutCC} corresponding to the potential $q(x)$. Notice that the matrix $\Phi(\lambda,x)$ is an entire function
of $\lambda$ by considering its Taylor series.

We saw in Proposition \ref{prop:SY}, that multiplication on a symmetric function $Y(\lambda)$ of the 
variable $\lambda$, which commute with $\Phi(\lambda,x,x_0)$ does not change vessel  parameters. In the case
of Srutm-Liouville vessel parameters we can describe these functions explicitly.
In order to understand which symmetric $x$-independent $Y(\lambda)$ commute with $\Phi(\lambda,x,x_0)$,
it is necessary and sufficient to understand when $Y(\lambda)$ commutes with its "generator" $\sigma_1^{-1}(\sigma_2 \lambda + \gamma)$. Extracting condition on $Y(\lambda)$ so that
\[ Y(\lambda) \sigma_1^{-1}(\sigma_2 \lambda + \gamma) = \sigma_1^{-1}(\sigma_2 \lambda + \gamma) Y(\lambda),
\] 
we will obtain that
\begin{equation} \label{eq:InCommPhi}
 Y(\lambda) = I + \bbmatrix{a(\lambda)& \dfrac{i c(\lambda)}{\lambda} \\ c(\lambda) & a(\lambda)}
\end{equation}
for functions $a(\lambda), c(\lambda)$, which are zero at infinity. It turns out that for a given $\gamma_*(x)$
in SL case, any two functions corresponding to 
the same vessel parameters differ by a constant symmetric function $Y(\lambda)$:
\begin{thm} \label{thm:UniquenessOfS}
Given SL vessel parameters $\sigma_1, \sigma_2, \gamma,\gamma_*(x)$, there exists a unique
initial value $S(0,\lambda)$ up to multiplication from the right on a symmetric, $x$-independent,
commuting with $\Phi(\lambda,x,x_0)$ function.
\end{thm}
\noindent\textbf{Proof:}
Given now two functions $S_1(\lambda,x)$, $S_2(\lambda,x)$ as in the theorem, the function 
$S_1^{-1}(\lambda,x)S_2(\lambda,x)$ will intertwine solutions of the input LDE with itself. Let us show that
such a function must be $x$-independent and commuting with $\Phi(\lambda,x,x_0)$. 

Using $\lambda = i s^2$ we find that
\[ \sigma_1^{-1} (\sigma_2 \lambda + \gamma) = \bbmatrix{0 & i \\ \lambda & 0},
\quad \Phi(\lambda,x,x_0) = V \bbmatrix{e^{-s (x-x_0)} & 0 \\ 0 & e^{s (x-x_0)}} V^{-1},
\]
where
\[ V = \bbmatrix{-\dfrac{1}{s} & \dfrac{1}{s} \\ 1 & 1}.
\]
Consequently, for the expression $S(\lambda,x) = \Phi(\lambda,x,x_0) S_0(\lambda) \Phi^{-1}(\lambda,x,x_0)$
to  be identity for $\lambda=\infty$, it is necessary to "cancel" the essential singularity arising from two
entire functions $\Phi(\lambda,x,x_0)$ and $\Phi^{-1}(\lambda,x,x_0)$ (or more precisely, let them cancel
each other). 

Using the formula for $\Phi(\lambda,x,x_0)$ in
\begin{eqnarray*}
\Phi(\lambda,x,x_0) S_0(\lambda,x_0) \Phi^{-1}(\lambda,x,x_0) = \\
V \bbmatrix{e^{-s (x-x_0)} & 0 \\ 0 & e^{k (x-x_0)}} V^{-1} S(\lambda,x_0)
V^{-1} \bbmatrix{e^{s (x-x_0)} & 0 \\ 0 & e^{-k (x-x_0)}} V^{-1} = I
\end{eqnarray*}
and considering coefficients of the exponents, we conclude that
\[ V^{-1} S(\lambda,x_0) V = \bbmatrix{b(s) & 0 \\ 0 & d(s)},
\]
for some analytic in $s$ at infinity functions $b(s), d(s)$. From here it follows that
\[ S(\lambda,x_0) = \bbmatrix{-\dfrac{b(s) + d(s)}{s} &  \dfrac{d(s) - b(s)}{s^2} \\ 
d(s) - b(s) & -\dfrac{b(s) + d(s)}{s} }.
\]
so, that $a(\lambda) = -\dfrac{b(s) + d(s)}{s}$ and
$c(\lambda) = d(s) - b(s)$ and we shall obtain that $S_0(\lambda)$ is of the form \eqref{eq:InCommPhi},
i.e. commutes with the fundamental matrix $\Phi(\lambda,x,x_0)$.
\qed

\subsection{Construction of $S(\lambda,x)$ for a given $\gamma_*(x)$. Classical case.}
Let us consider the Sturm-Liouville differential equation
\eqref{eq:SL}
\[ -y''(x,s) + q(x) y(x,s) = s^2 y(x,s),
\]
where the potential $q(x)$ is a real measurable function satisfying the condition
\cite{bib:FaddeyevII}
\[ \int\limits_{-\infty}^\infty (1 + |x|) |q(x)| dx < \infty.
\]
Consider the following solutions $f_1(x,s), f_2(x,s)$ defined from a Volterra type
equation ($\Im s > 0$)
\begin{eqnarray}
\label{eq:f1} f_1(x,s) = e^{ikx} + \int\limits_{-\infty}^\infty G_1(x-y,s) q(y) f_1(y,s) dy, \\
\label{eq:f2} f_2(x,s) = e^{-ikx} + \int\limits_{-\infty}^\infty G_2(x-y,s) q(y) f_2(y,s) dy,
\end{eqnarray}
where
\[ G_1(t,s) = - Hev(-t) \dfrac{\sin(st)}{s}, \quad G_2(t,s) = Hev(t) \dfrac{\sin(s t)}{s}
\]
and $Hev(t)$ is the Heaviside function
\[ Hev(t)=1, t>0, \quad Hev(t)=0, t<1.
\]
The functions $f_1, f_2$ behave \cite[1.4, 1.5]{bib:FaddeyevII} as $e^{isx}$ 
and $e^{-isx}$ for $x$ approaching $+\infty$ and $-\infty$ respectively. Moreover, the following bounds hold
\cite{bib:Levinson}
\begin{eqnarray}
\label{eq:f1Bound} |f_1(x,s) - e^{isx}| \leq 
C \dfrac{e^{-\Im s x}}{1 + |s|} \int\limits_x^\infty ( 1 + |y|) |q(y)| dy, \\
\label{eq:f2Bound} |f_1(x,s) - e^{-isx}| \leq 
C \dfrac{e^{\Im s x}}{1 + |s|} \int\limits_{-\infty}^x ( 1 + |y|) |q(y)| dy.
\end{eqnarray}
Then Z. S. Agranovich and V. A. Marchenko \cite{bib:AgMarch} proved that the same solutions $f_1, f_2$
may be represented as ($\Im s > 0$)
\[ f_1(x,s) = e^{isx} + \int\limits_x^\infty A_1(x,y) e^{isy} dy, \quad
f_2(x,s) = e^{-isx} + \int\limits_x^\infty A_2(x,y) e^{-isy} dy
\]
where $A_1, A_2$ are square integrable functions of $y$ for each $x$. Moreover they also showed that defining
\[ \xi_1(x) = \int\limits_{x}^\infty |q(y)| dy, \quad
\xi_2(x) = \int\limits_{-\infty}^{x} |q(y)| dy,
\]
the functions $A_1, A_2$ satisfy Volterra type equations \cite[1.10, 1.11]{bib:FaddeyevII} 
and successive approximations give the following bounds for them \cite[1.12]{bib:FaddeyevII} 
\[ |A_1(x,y)| \leq C \xi_1 (\dfrac{x+y}{2}), \quad
|A_2(x,y)| \leq C \xi_2 (\dfrac{x+y}{2}).
\]
One can also find that 
\begin{eqnarray}
\label{eq:A1Bound}|\dfrac{\partial}{\partial x} A_1(x,y) + \dfrac{1}{4} q(\dfrac{x+y}{2})| \leq 
C \xi_1(x) \xi_1(\dfrac{x+y}{2}), \\
\label{eq:A2Bound}|\dfrac{\partial}{\partial x} A_2(x,y) - \dfrac{1}{4} q(\dfrac{x+y}{2})| \leq 
C \xi_2(x) \xi_2(\dfrac{x+y}{2}),
\end{eqnarray}
and
\[ q(x) = - 2 \dfrac{\partial}{\partial x} A_1(x,x) = 
2 \dfrac{\partial}{\partial x} A_2(x,x).
\]
Finally, we define
\begin{equation} \label{eq:DefBeta}
\beta(x) = - A_1(x,x) = - \dfrac{1}{2} \int\limits_x^\infty q(y) dy
\end{equation}
and
\[ S(\lambda,x) = \Phi_*(x,\lambda) \Phi^{-1}(x,\lambda),
\]
where $\Phi, \Phi_*$ are solutions of the input and the output LDE respectively,
corresponding to the SL parameters, defined using the function $\beta(x)$. One can take
\[ \Phi(x,\lambda) = \bbmatrix{ \cos(s x)& \dfrac{i \sin(sx)}{s}\\ i s \sin(sx) & \cos(sx)}
\]
and
\[ \Phi_*(x,\lambda) = \bbmatrix{\dfrac{f_1 + f_2}{2}& \dfrac{f_1 - f_2}{2s}\\
\dfrac{f_1' + f_2'}{2i}-\beta\dfrac{f_1 + f_2}{2i} & \dfrac{f_1' - f_2'}{2is}-\beta\dfrac{f_1 - f_2}{2is}}.
\]
Note that the function $\Phi$ depend on $s^2=-i\lambda$ and as a result is an entire function of $\lambda$.
Simple calculations show that
\begin{multline*}
S(\lambda,x)
= \bbmatrix{\dfrac{f_1e^{-isx} + f_2e^{isx}}{2}& \dfrac{f_1e^{-isx} - f_2e^{isx}}{2s}\\
\dfrac{(f_1'-\beta f_1)e^{-isx} + (f_2' - \beta f_2)e^{isx}}{2i} & 
\dfrac{(f_1'-\beta f_1)e^{-isx} - (f_2' - \beta f_2)e^{isx}}{2is}}
\end{multline*}
\begin{thm} For all $\lambda=-is^2$ with $\Im s>\epsilon>0$ it holds that 
\[ \lim\limits_{\lambda \rightarrow\infty} S(\lambda,x) = I.
\]
\end{thm}
\textbf{Proof:} For the first row of $S$ notice that from \eqref{eq:f1Bound}, \eqref{eq:f2Bound} it follows that
\[ f_1(x,s) e^{-isx} = 1 + o(1), \quad
f_2(x,s) e^{isx} = 1 + o(1),
\]
where $o(1)$ means a function going to zero as $s$ goes to infinity for all $x$. For the second row, let us consentrate
first on the derivative of $f_1(x,s)$, which can be found from \eqref{eq:f1}
\[ f_1'(x,s) = is e^{isx} - A_1(x,x) e^{isx} + 
\int\limits_{x}^\infty \dfrac{\partial}{\partial x} A_1(x,y) e^{isy} dy.
\]
From where it follows that
\[ f_1'(x,s) e^{-isx} - is  + A_1(x,x) =
\int\limits_{x}^\infty \dfrac{\partial}{\partial x} A_1(x,y) e^{is(y-x)} dy
\]
Then one can estimate using \eqref{eq:A1Bound} that for $\Im s > 0$
\[ \begin{array}{llll}
|\int\limits_{x}^\infty \dfrac{\partial}{\partial x} A_1(x,y) e^{is(y-x)} dy| \leq \\
\leq \int\limits_{x}^\infty [ |\dfrac{1}{4} q(\dfrac{x+y}{2})| +
C \xi_1(x) \xi_1(\dfrac{x+y}{2})] e^{-\Im s(y-x)} dy  \\
\leq \int\limits_{x}^\infty [ |\dfrac{1}{4} q(\dfrac{x+y}{2})| dy \int\limits_{x}^\infty  e^{-\Im s(y-x)} dy
+ C \xi_1(x) \xi_1(-\infty) \int\limits_{x}^\infty e^{-\Im s(y-x)} dy  \\
\leq K(x) \int\limits_{x}^\infty e^{-\Im s(y-x)} = K(x) \dfrac{1}{\Im s},
\end{array} \]
from where it follows that when $\Im s \rightarrow +\infty$ the integral
$\int\limits_{x}^\infty \dfrac{\partial}{\partial x} A_1(x,y) e^{is(y-x)} dy$ approaches zero.
For the case when $\Im s>0$ is fixed, we notice that then 
\[ \int\limits_{x}^\infty \dfrac{\partial}{\partial x} A_1(x,y) e^{is(y-x)} dy =
\int\limits_{x}^\infty \dfrac{\partial}{\partial x} A_1(x,y) e^{-\Im s(y-x)} 
e^{i \Re s (y-x)} dy
\]
and from the above calculation it follows that $\dfrac{\partial}{\partial x} A_1(x,y) e^{-\Im s(y-x)}$
is $L^1(\mathbb R)$ function. By the Riemann-Lebesgue lemma as $\Re s \rightarrow \pm \infty$, the integral
approaches zero too. Thus we conclude that
\[ f_1'(x,s) e^{-isx} - is  + A_1(x,x) = o(1),
\]
and as a result, by the definition \eqref{eq:DefBeta} of $\beta$
\[ (f_1' - \beta f_1) e^{-isx} = is  - A_1(x,x) - \beta f_1 e^{-isx} + o(1) = is + o(1).
\]
Similarly one can show that
\[ (f_2' - \beta f_2) e^{isx} = -is  + A_2(x,x) - \beta f_2 e^{2isx} + o(1) = -is + o(1).
\]
From where it follows the statement for the second row of $S(\lambda,x)$.
\qed

Finally, we focus on the function $S(\lambda,x)$ at the value $x=0$. Then notice that
$f_i(0,s) = \bar f_i(0,-\bar s)$ ($i=1,2$). As a result we can define a function
\[ M_i(\mu) = f_i(0,\sqrt{\mu}),
\]
where we choose the root in such a manner that $\Im \sqrt{\mu} = s \geq 0$.
Consequently, for $\mu = s^2$ it holds that $\sqrt{\bar \mu} = -\bar s$, which must be at
the upper half plane, and
\[ M_i(\bar \mu) = f_i (0,\sqrt{\bar \mu}) = f_i (0,-\bar s) = \bar f_i(0,s)
= \bar M_i(\mu).
\]
Consequently, the function $M_i(\lambda)$ is bounded for all $\lambda$, has the value $1$ at infinity
(for big $\mu$), analytic at the whole complex plane except for a cut on the 
positive real axis, where it has a jump. Similarly, we define the functions, corresponding to the derivatives
of $f_i(x,s)$, $M_i^1(0,\mu) = f_i'(0,\sqrt{\mu})$.

Finally, we substitute $\mu = -i\lambda$ and define the function 
\[ S(\lambda,0) = \bbmatrix{\dfrac{M_1(-i\lambda) + M_2(-i\lambda)}{2} &
\dfrac{M_1(-i\lambda) - M_2(-i\lambda)}{2\sqrt{-i\lambda}} \\
\dfrac{M_1^1(-i\lambda) + M_2^1(-i\lambda)}{2} &
\dfrac{M_1^1(-i\lambda) - M_2^1(-i\lambda)}{2\sqrt{-i\lambda}}
}
\]
where again, the function $\sqrt{-i\lambda}$ may be defined except for a cut on the imaginary axis.
Notice that
\[ S(\lambda,x) = \Phi_*(x,\lambda) S(\lambda,0) \Phi^{-1}(\lambda,x),
\]
Where $\Phi_*(x,\lambda)$ is the fundamental solution of the output LDE, attaining $I$ at $x=0$.

Let us check next the symmetric condition. Notice that from Theorem \ref{thm:S=S-1*} it is enough to check the symmetry
of $S(\lambda,0)$ only. In order to do it, we notice that all the entries of $S(\lambda,0)$ are created
from the functions $f_i(x,s)$ satisfying $\bar f_i(x,-\bar s) = f_i(x,s)$. As a result, it holds that
\[ \bar M_i(i\bar \lambda) = \bar M_i(\bar \mu) = M_i(\mu) = M_i(-i\lambda)
\]
and similarly, $\bar M_i^1(i\bar\lambda) = M^1_i(-i\lambda)$. Thus denoting
$S(\lambda,0) = \bbmatrix{a(\lambda)&b(\lambda)\\c(\lambda) & d(\lambda)}$
we obtain that
\[ \bar a(-\bar\lambda) = a(\lambda), \quad \bar b(-\bar\lambda) = -b(\lambda)
\]
and
\[ \bar d(-\bar\lambda) = d(\lambda), \quad \bar c(-\bar\lambda) = -c(\lambda).
\]
As a result the symmetry condition for $S(\lambda,0)$ becomes
\[ \begin{array}{llllll}
S^*(-\bar\lambda,0) \sigma_1 S(\lambda,0) =
\bbmatrix{\bar a(-\bar\lambda) & \bar c(-\bar\lambda)\\\bar b(-\bar\lambda) & \bar d(-\bar\lambda)}
\bbmatrix{0&1\\1&0} \bbmatrix{a(\lambda) & b(\lambda)\\c(\lambda)&d(\lambda)} \\
= \bbmatrix{\bar c(-\bar\lambda) a(\lambda) + \bar a(-\bar\lambda)c(\lambda) &
 \bar c(-\bar\lambda) b(\lambda) + \bar a(-\bar\lambda)d(\lambda) \\
 \bar d(-\bar\lambda) a(\lambda) + \bar b(-\bar\lambda)c(\lambda) &
 \bar d(-\bar\lambda) b(\lambda) + \bar b(-\bar\lambda)d(\lambda)} \\
= \bbmatrix{-c(\lambda) a(\lambda) + a(\lambda)c(\lambda) &
 -c(\lambda) b(\lambda) + a(\lambda)d(\lambda) \\
 d(\lambda) a(\lambda) - b(\lambda)c(\lambda) &
 d(\lambda) b(\lambda) - b(\lambda)d(\lambda)} \\
 = \bbmatrix{0 & -c(\lambda) b(\lambda) + a(\lambda)d(\lambda) \\ d(\lambda) a(\lambda) - b(\lambda)c(\lambda) & 0} \\
= [d(\lambda) a(\lambda) - b(\lambda)c(\lambda)] \bbmatrix{0&1\\1&0}
\end{array} \]
Thus we obtain
\begin{thm} The function
\[ \dfrac{1}{\det S(\lambda,0)} S(\lambda,x)
\]
is in the class $\Ir(\sigma_1, \sigma_2, \gamma, \gamma_*(x), \mathrm I =\mathbb R)$. In other words,
it realizes the given $\gamma_*(x)$.
\end{thm}

The following theorem was established for the finite dimensional case in \cite{bib:SLVessels}. We present
now its generalization
\begin{thm} \label{prop:Gamma*Formula}
For Sturm Liouville vessel the following formula for $\gamma_*$ holds
\[ \gamma_* = \gamma + \bbmatrix{i\frac{\tau''}{\tau} & \frac{\tau'}{\tau} \\ -\frac{\tau'}{\tau} & 0}
\]
\end{thm}
\noindent\textbf{Proof:} Let us take 
\[ \gamma_* = \bbmatrix{-i (\beta'(x) - \beta^2(x)) & -\beta(x) \\ \beta(x) & i}
\]
If we denote $-\beta=\dfrac{\tau'}{\tau}$, then $\gamma_* = \bbmatrix{i \dfrac{\tau''}{\tau} & \dfrac{\tau'}{\tau} \\ -\dfrac{\tau'}{\tau} & i}$ and we have to prove that $-\beta = \dfrac{\tau'}{\tau}$. Consider now the formula
\eqref{eq:GKform}
\begin{equation} \label{eq:TauBeta}
 \dfrac{\tau'}{\tau}  = \TR (\sigma_2 B^*(x) \mathbb X^{-1}(x)B(x)).
\end{equation}
Notice that the expression $B^*(x) \mathbb X^{-1}(x)B(x)$ is the first moment and from the general
formula of the first moment appearing in \cite[section 5]{bib:SchurVessels} (in this article the moments are
different from the moments in this article by multiplication on $\sigma_1^{-1}$ from the right)
\[ H_0(x) = \bbmatrix{-\beta & \dfrac{r + i(\beta'-\beta^2)}{2} \\ \dfrac{r - i(\beta'+\beta^2)}{2} & h_0^{21} }
\]
it follows that
\begin{equation} \label{eq:BetaTau}
 \dfrac{\tau'}{\tau} = \TR (\sigma_2 H_0(x)) = -\beta
\end{equation}
as desired.
\qed
\begin{cor} The following formula for the potential holds
\begin{equation} \label{eq:qtau}
q(x) = - 2 \dfrac{d^2}{dx^2} \ln \tau(x).
\end{equation}
\end{cor}
\noindent\textbf{Proof:} Immediate from \eqref{eq:qbeta}.
\qed

This notion is important in the sense that conjecturally the matrix $S(\lambda,x), q(x)$ and the solutions of
\eqref{eq:OutCC} may be represented from it and $e^{isx}$, which are solutions of the input LDE \eqref{eq:InCC}.
We make the following conjecture (generalization of \cite[theorem 3]{bib:SLVessels})
\begin{conj} The entries of the matrix $S(\lambda,x)$ are linear combinations of $\dfrac{\tau^{(n)}(x)}{\tau(x)}$
in some $p$-norm on $\mathrm I$, in other words, the entries are of the form
\[ \sum \alpha_n \dfrac{\tau^{(n)}(x)}{\tau(x)}, \quad \alpha_n\in\mathbb C.
\]
\end{conj}

\subsection{Scattering data versus $S(\lambda,x_0)$. Gelfand-Levitan equation.}
Following \cite{bib:FadeevInv} for the case $\int\limits_0^\infty x |q(x)| dx<\infty$ there are introduced Jost solutions
$\phi(x,s)$ and $f(x,s)$ \cite{bib:Jost, bib:Levinson}
\begin{eqnarray}
\label{eq:phiCond}	\phi(x,s) & : & \phi(0,s) =0, \quad \phi'(0,s) = 1, \\
\label{eq:fCond}		f(x,s) & : & \lim\limits_{x\rightarrow\infty} e^{-isx} f(x,s) = 1.
\end{eqnarray}
and defining further $M(k) = \phi'(x,k) f(x,s) - f'(x,s) \phi(x,k)$ one reconstructs the potential $q(x)$
using Gelfand-Levitan equation \cite{bib:GL} (or alternatively Marchenko equation \cite{bib:Marchenko}). There are two steps, essential for
this construction, namely, one considers the case when the spectrum of the operator $L$ is purely continuous
and the case when this spectrum additionally contains finite number of points.
For the purely continuous case, one proves that there is a solution $K(x,y) $ of
the Gelfand-Levitan equation
\cite[(8.5)]{bib:FadeevInv}
\begin{equation} \label{eq:GelfandLevitan}
K(x,y) + \Omega(x,y) + \int\limits_0^x K(x,t) \Omega(t,y) dt = 0, \quad x>y.
\end{equation}
where $\Omega(x,y)$ is uniquely defined from $M(k)$ by \cite[(8.4)]{bib:FadeevInv}
\[ \Omega(x,y) = 2/\pi \int_0^\infty \dfrac{\sin(kx)}{k}[\dfrac{1}{M(k)M(-k)}-1] \dfrac{\sin(ky)}{k} k^2 dk.
\]
The formula for the potential is \cite[(10.4)]{bib:FadeevInv}
$q(x) = 2 \dfrac{d}{dx} K(x,x)$. Then one make a modification, so that the discrete spectrum is taken into
account \cite[(8.14, 8.15)]{bib:FadeevInv}. We are going to present analogues of these formulas in our setting.
Suppose that we are given a vessel \eqref{eq:DefKV}
\[ \mathfrak{K_V} = (A, B(x), \mathbb X(x); \sigma_1(x), 
\sigma_2(x), \gamma(x), \gamma_*(x);
\mathcal{H}, \mathcal{E};\mathrm I), 
\]
and let as fix an arbitrary $x_0\in\mathrm I$. Define
\begin{equation} \label{eq:OmeagaDef}
\Omega(x,y) = \bbmatrix{1&0} B^*(x) \mathbb X^{-1}(x_0) B(y) \bbmatrix{1\\0}.
\end{equation}
and
\begin{equation} \label{eq:KDef}
K(x,y) = -\bbmatrix{1&0} B^*(x) \mathbb X^{-1}(x) B(y) \bbmatrix{1\\0}.
\end{equation}
Then Gel'fand-Levitan equation \eqref{eq:GelfandLevitan} holds
\[ \begin{array}{llllll}
K(x,y) + \Omega(x,y) + \int\limits_{x_0}^x K(x,t) \Omega(t,y) dt = \\
= K(x,y) + \Omega(x,y) - \int\limits_{x_0}^x \bbmatrix{1&0} B^*(x) \mathbb X^{-1}(x) B(t) \bbmatrix{1\\0}  \bbmatrix{1&0} B^*(t) \mathbb X^{-1}(x_0) B(y) \bbmatrix{1\\0}dt = \\
= K(x,y) + \Omega(x,y) - \bbmatrix{1&0} B^*(x) \mathbb X^{-1}(x)  \int\limits_{x_0}^x B(t) \sigma_2 B^*(t) dt \mathbb X^{-1}(x_0) B(y) \bbmatrix{1\\0} = \\
= \text{ using vessel condition \eqref{eq:DB} } = \\
= K(x,y) + \Omega(x,y) - \bbmatrix{1&0} B^*(x) \mathbb X^{-1}(x) (\mathbb X(x) - \mathbb X(x_0)) \mathbb X^{-1}(x_0) B(y) \bbmatrix{1\\0} = \\
= K(x,y) + \Omega(x,y) - \bbmatrix{1&0} B^*(x) \mathbb X^{-1}(x_0) B(y) \bbmatrix{1\\0} +
\bbmatrix{1&0} B^*(x) \mathbb X^{-1}(x)  B(y) \bbmatrix{1\\0} = 0.
\end{array} \]
Finally, the formula \eqref{eq:qbeta} for the potential gives
\[ \begin{array}{llllll}
q(x) = 2\dfrac{d}{dx} \SP ( \mathbb X^{-1}(x)  \dfrac{d}{dx} \mathbb X(x)) = 
\SP ( \mathbb X^{-1}(x)  B(x) \sigma_2 B^*(x)) = \\
\quad = 2\dfrac{d}{dx} \SP ( \mathbb X^{-1}(x)  B(x) \bbmatrix{1\\0} \bbmatrix{1&0} B^*(x)) = 
\dfrac{d}{dx}\big( \bbmatrix{1&0} B^*(x) \mathbb X^{-1}(x)  B(x) \bbmatrix{1\\0} \big) = \\
\quad = 2\dfrac{d}{dx} K(x,x),
\end{array} \]
which is identical to \cite[(10.4)]{bib:FadeevInv}.

\subsection{\label{sec:GlobalScat}Jost solutions.}
In the sequel we will use the following number
\begin{equation} \label{eq:mA}
 m(A) = \max \{ \Im \lambda \mid \lambda\in\SPEC(A)\}.
\end{equation}
\begin{lemma} \label{lemma:PosSig2}
For all $x\in\mathrm I$ it holds that $\mathbb X(x) \geq \mathbb X_0$. The operator $B(x)$ satisfies the following
\[ \lim\limits_{x\rightarrow\infty} B(x) e^{i s x} = 0, \quad \Im s> m(A) .
\]
\end{lemma}
\noindent\textbf{Proof:} The inequality is immediate from the formula \eqref{eq:DX}
\[ \mathbb X(x) = \mathbb X_0 + \int\limits_{x_0}^x B(y)\sigma_2 B(y) dy,
\]
since $\sigma_2\geq 0$ and positive operators form a convex set inside the space of all operators. 
For the second part, we can use Dunford-Schwartz calculus \cite{bib:DanSchw}
\[ \begin{array}{lll}
B(x) e^{i s x} = \dfrac{1}{2\pi i} \oint (\lambda I - A)^{-1} B_0 \Phi^t (\lambda,x,x_0) d\lambda e^{isx} = \\
= \dfrac{1}{2\pi i} \oint (\lambda I - A)^{-1} B_0 
\bbmatrix{\cos(t x) & i t \sin(t x)\\ \dfrac{i \sin(t x)}{t} & \cos(t x)} d\lambda e^{isx} , \quad \lambda=it^2.
\end{array} \]
Taking the norm of this expression, we shall obtain that
\[ \|B(x) e^{i s x}\| \leq
\dfrac{1}{2\pi i}  \oint \| (\lambda I - A)^{-1} B_0 \| \max_{t} |e^{i(s\pm t)x}| d\lambda \leq C(s) 
e^{-(\Im s - m(A)) x},
\]
where $C(s)$ is a constant function, depending on $s$ only. When $x$ tend to infinity, we obtain the desired result.
\qed

Using the transfer function $S(x,s)$, we will look for the Jost solutions,
defined from the following formulas
\begin{eqnarray}
\label{eq:phiDef}\phi(x,s) & = \bbmatrix{1 & 0} S(x,s) \Phi(x,s) \bbmatrix{\phi_1(s)\\\phi_2(s)}, \\
\label{eq:fDef}f(x,s) & = \bbmatrix{1&0}S(x,s) \Phi(x,s) \bbmatrix{1 \\ s} f(s),
\end{eqnarray}
where $\phi_1(s), \phi_2(s), f(s)$ are functions, which must be found.
Solving for the $\phi_1, \phi$ so that the condition \eqref{eq:phiDef} are fulfilled,
one will come to the conclusion that
\[ \phi(x,s) = \bbmatrix{1 & 0} \Phi_*(\lambda,x,x_0) \bbmatrix{0\\-i}.
\]
Particularly, the Jost solution $\phi(\lambda,x)$ is an entire function of $\lambda$ for each $x$. 
The choice for the function $f(x,s)$ comes from the following formula:
\[ f(x,s) =
\bbmatrix{1&0} S(x,s) \Phi(\lambda,x,x_0) \bbmatrix{1 \\ s} f(s) =
\bbmatrix{1&0}S(x,s)  \bbmatrix{e^{isx} \\ s e^{isx}} f(s).
\]
So if we want to satisfy the condition of the Jost solution, we will demand that
\[  \lim\limits_{x\rightarrow\infty} e^{-isx} f(x,s) = 
\lim\limits_{x\rightarrow\infty} 
\bbmatrix{1&0}S(x,s)  \bbmatrix{1 \\ s } f(s) = 1.
\]
Let us define $f_1(x,s) = \bbmatrix{1&0} S(x,s) \Phi(\lambda,x,x_0) \bbmatrix{1 \\ s}$, then the following
lemma holds
\begin{lemma} Suppose that we are given a vessel $\mathfrak{K_V}$ on a half axis 
$\mathrm I =[x_0,\infty)$, then for $\Im s > m(A)$ it holds that
\[ \int\limits_0^\infty |f_1(y,s)|^2 dy = -\dfrac{\bbmatrix{1&\bar s} S^*(\lambda,0)\sigma_1 S(\lambda,x) \bbmatrix{1\\s}}{\lambda + \bar\lambda}.
\]
\end{lemma}
\textbf{Proof:}  We shall use the fundamental matrices
$\Phi_* = \Phi_*(s,x,x_0)$, and $\Phi = \Phi(s,x,x_0)$ of the output and the input LDEs respectively. Then
\[ \begin{array}{llll}
f_1^*(x,s) f_1(x,s) & = \bbmatrix{1 & \bar s} \Phi^* S^*(\lambda,x)\bbmatrix{1 \\ 0}\bbmatrix{1&0} S(\lambda,x) \Phi \bbmatrix{1 \\ s} =\\
& = \bbmatrix{1 & \bar s}  \Phi^* S^*(\lambda,x)\sigma_2 S(\lambda,x)\Phi \bbmatrix{1 \\ s} =\\
& = \bbmatrix{1 & \bar s}  S^*(\lambda,0) \Phi_*^* \sigma_2 \Phi_* S(\lambda,x) \bbmatrix{1 \\ s} =\\
& = \bbmatrix{1 & \bar s}  S^*(\lambda,0) \dfrac{\partial}{\partial x} \big( \dfrac{\Phi_*^* \sigma_1 \Phi_*}{\lambda+\bar\lambda} \big) S(\lambda,0) \bbmatrix{1 \\ s}.
\end{array} \]
integrating the last equation and using additionally \eqref{eq:SInttw} we obtain that
\[ \begin{array}{llll}
\int\limits_0^x |f_1(y,s)|^2 dy & = 
\bbmatrix{1 & \bar s}  S^*(\lambda,0) [\dfrac{\Phi_*^* \sigma_1 \Phi_*}{\lambda+\bar\lambda} -
\dfrac{\sigma_1}{\lambda+\bar\lambda}] S(\lambda,0) \bbmatrix{1 \\ s} = \\
& = \bbmatrix{1 & \bar s}  \dfrac{\Phi^* S^*(\lambda,x) \sigma_1 S(\lambda,x)\Phi}{\lambda+\bar\lambda} \bbmatrix{1 \\ s} - \bbmatrix{1 & \bar s}  \dfrac{ S^*(\lambda,0) \sigma_1 S(\lambda,0)}{\lambda+\bar\lambda} \bbmatrix{1 \\ s}.
\end{array} \]
Using here the expression \eqref{eq:Phi} for $\Phi$:
\begin{equation} \label{eq:f1Before}
 \begin{array}{llll}
\int\limits_0^x |f_1(y,s)|^2 dy & = 
\bbmatrix{1 & \bar s}  \dfrac{S^*(\lambda,x) \sigma_1 S(\lambda,x)}{\lambda+\bar\lambda} \bbmatrix{1 \\ s}e^{i(s-\bar s)x} - \bbmatrix{1 & \bar s}  \dfrac{ S^*(\lambda,0) \sigma_1 S(\lambda,0)}{\lambda+\bar\lambda} \bbmatrix{1 \\ s} \\
& = \bbmatrix{1 & \bar s}  \dfrac{S^*(\lambda,x) \sigma_1 S(\lambda,x) - \sigma_1}{\lambda+\bar\lambda} \bbmatrix{1 \\ s}e^{i(s-\bar s)x} + \\
& \hspace{3cm} + \bbmatrix{1 & \bar s}  \dfrac{\sigma_1e^{i(s-\bar s)x} - S^*(\lambda,0) \sigma_1 S(\lambda,0)}{\lambda+\bar\lambda} \bbmatrix{1 \\ s}.
\end{array} \end{equation}
Let us focus on the expression 
\[ \bbmatrix{1 & \bar s}  \dfrac{\sigma_1 - S^*(\lambda,x) \sigma_1 S(\lambda,x)}{\lambda+\bar\lambda} \bbmatrix{1 \\ s} = \sigma_1 B^*(x) (\bar\lambda I - A^*)^{-1} \mathbb X^{-1}(x) (\lambda I - A)^{-1} B(x) \sigma_1.
\] 
From Lemma \ref{lemma:PosSig2}, for all $x\in\mathrm I$ it holds that $\mathbb X(x) \geq \mathbb X(x_0)$ so
\[ \sigma_1 B^*(x) (\bar\lambda I - A^*)^{-1} \mathbb X^{-1}(x) (\lambda I - A)^{-1} B(x) \sigma_1 \leq
\sigma_1 B^*(x) (\bar\lambda I - A^*)^{-1} \mathbb X^{-1}(x_0) (\lambda I - A)^{-1} B(x) \sigma_1
\]
Taking the norm, we shall obtain that
\[ \begin{array}{lll}
\| \sigma_1 B^*(x) (\bar\lambda I - A^*)^{-1} \mathbb X^{-1}(x) (\lambda I - A)^{-1} B(x) \sigma_1 \| 
& \leq \| \sigma_1 B^*(x) (\bar\lambda I - A^*)^{-1} \mathbb X^{-1}(x_0) (\lambda I - A)^{-1} B(x) \sigma_1 \| \\
& \leq K(\lambda) \| B(x) \|^2, \quad \text{$K(\lambda)$ - x - independent}. 
\end{array} \]
Multiplying this inequality by $e^{i(s-\bar s)}$ and using condition on $\|B(x)\|$ in Lemma \ref{lemma:PosSig2}, we
obtain that for $\Im s > m(A)$
\[ \lim\limits_{x\rightarrow\infty}
\| \bbmatrix{1 & \bar s}  \dfrac{\sigma_1 - S^*(\lambda,x) \sigma_1 S(\lambda,x)}{\lambda+\bar\lambda} \bbmatrix{1 \\ s} e^{i(s-\bar s)x}\| \leq 0
\]
Plugging this into equation \eqref{eq:f1Before} we obtain that
\[ \int\limits_0^\infty |f_1(y,s)|^2 dy = - \bbmatrix{1 & \bar s}  \dfrac{ S^*(\lambda,0) \sigma_1 S(\lambda,0)}{\lambda+\bar\lambda} \bbmatrix{1 \\ s}.
\]
\qed

Let us define the following expression, which is essential for the existence of the Jost solution $f(x,s)$
\begin{equation} \label{eq:h}
h(x,s) = \dfrac{f(x,s)}{f(s)} e^{-isx} = \bbmatrix{1&0}S(x,s)  \bbmatrix{1 \\ s }
\end{equation}
Let as also denote
\begin{equation} \label{eq:KS}
 K_S(x,s) = \dfrac{\bbmatrix{1 & \bar s}S^*(\lambda,x) \sigma_1 S(\lambda,x)\bbmatrix{1 \\ s} }{\lambda+\bar\lambda},
\end{equation}
then the following theorem holds:
\begin{thm} \label{thm:hisOK} Suppose that we are given a vessel
\eqref{eq:DefKV}
\[ \mathfrak{H_V} = (A, B(x), \mathbb X(x); \sigma_1(x), 
\sigma_2(x), \gamma(x), \gamma_*(x);
\mathcal{H}, \mathcal{E};\mathrm I = (x_0,\infty)), 
\] 
then the function $h(x,s) = |h(x,s)| e^{i\Theta_h(x,s)}$ defined in \eqref{eq:h} has the following properties
\begin{enumerate}
	\item $h^*(x,-\bar s) = h(x,s) \det S(\lambda,x_0)$,
	\item $|h(x,s)|^2 = \dfrac{\partial}{\partial x} K_S(s,x) + i(s-\bar s)K_S(x)$,
	\item $2 \dfrac{\partial}{\partial x} \Theta_h(x,s) = -\dfrac{s+\bar s}{|h(x,s)|^2} \dfrac{\partial}{\partial x} K_S(x,s)$.
\end{enumerate}
\end{thm}
\noindent\textbf{Proof: 1.} Using the definition of $h(x,s)$ and Lemma \ref{lemma:detS} we obtain that
\[ \begin{array}{llllll}
\bar h (x,-\bar s) & = \bbmatrix{1 & -s} S^*(x,-\bar\lambda) \bbmatrix{1\\0}
 = \bbmatrix{1 & -s} \sigma_1 S^{-1}(x,\lambda) \sigma_1^{-1} \bbmatrix{1\\0}  = \\
& = \bbmatrix{-s & 1} S^{-1}(x,\lambda) \bbmatrix{0\\1} =  \text{ denoting } S(x,\lambda) = \bbmatrix{a&b\\c&d}\\
& = \bbmatrix{-s & 1} \bbmatrix{d & -b \\ -c & a} \dfrac{1}{\det S(\lambda,x)} \bbmatrix{0\\1}  = \\
& = (a + s b) \dfrac{1}{\det S(\lambda,x)} = h(x,s) \dfrac{1}{\det S(\lambda,x)} = h(x,s) \dfrac{1}{\det S(\lambda,x_0)}.
\end{array} \]

\noindent\textbf{2.} Let us denote $\Phi=\Phi(\lambda,x,x_0), \Phi_*=\Phi_*(\lambda,x,x_0)$, then
\[ \begin{array}{llll}
h^*(x,s) h(x,s) & = \bbmatrix{1 & \bar s} S^*(\lambda,x)\bbmatrix{1 \\ 0}\bbmatrix{1&0} S(\lambda,x) \bbmatrix{1 \\ s} = \\
& =  \bbmatrix{1 & \bar s} S^*(\lambda,x)\sigma_2 S(\lambda,x) \bbmatrix{1 \\ s} = \\
& = \bbmatrix{1 & \bar s} \Phi^{-1*} S^*(\lambda,x_0) \Phi_*^* \sigma_2 
\Phi_* S(\lambda,x_0) \Phi^{-1} \bbmatrix{1 \\ s} = \text{using \eqref{eq:DPhi^*S2Phi}} \\
& = \bbmatrix{1 & \bar s} \Phi^{-1*} S^*(\lambda,x_0) \dfrac{\partial}{\partial x} \big( \dfrac{\Phi_*^* \sigma_1 \Phi_*}{\lambda+\bar\lambda} \big) S(\lambda,x_0) \Phi^{-1} \bbmatrix{1 \\ s}. 
\end{array} \]
Using the formula \eqref{eq:Phi} for $\Phi$, we can calculate that
$\Phi \bbmatrix{1 \\ s} = \bbmatrix{1 \\ s} e^{isx}$. As a result, using this and \eqref{eq:SInttw},
the formula for $h^*(x,s) h(x,s)= |h(x,s)|^2$ becomes
\[ \begin{array}{llll}
|h(x,s)|^2  & = \dfrac{\partial}{\partial x} \big(\dfrac{\bbmatrix{1 & \bar s}S^*(\lambda,x) \sigma_1 S(\lambda,x)\bbmatrix{1 \\ s} }{\lambda+\bar\lambda} e^{i (s-\bar s)x} \big)  e^{-i (s-\bar s)} = \\
& = \dfrac{\partial}{\partial x} \big(\dfrac{\bbmatrix{1 & \bar s}S^*(\lambda,x) \sigma_1 S(\lambda,x)\bbmatrix{1 \\ s} }{\lambda+\bar\lambda} \big) + i (s-\bar s)
\dfrac{\bbmatrix{1 & \bar s}S^*(\lambda,x) \sigma_1 S(\lambda,x)\bbmatrix{1 \\ s} }{\lambda+\bar\lambda} \\
& = \dfrac{\partial}{\partial x} K_S(s,x) + i(s-\bar s)K_S(x).
\end{array} \]

\noindent\textbf{3.} Using the formula \eqref{eq:DS} we find that
\[ \begin{array}{lll}
\bar h(x,s) h'(x,s) - \bar h'(x,s) h(x,s) =
\bbmatrix{1&\bar s} S^*(\lambda,x) \bbmatrix{1\\0} 
\bbmatrix{\beta(x)&i} S(\lambda,x) \bbmatrix{1\\s} - \\
\quad \quad - \bbmatrix{1&\bar s} S^*(\lambda,x) \bbmatrix{\beta(x)\\-i} 
\bbmatrix{1&0} S(\lambda,x) \bbmatrix{1\\s} - i(s+\bar s)  \bar h(x,s) h(x,s) = \\
\quad = i \bbmatrix{1&\bar s} S^*(\lambda,x) \sigma_1  S(\lambda,x) \bbmatrix{1\\s} - i(s + \bar s)  \bar h(x,s) h(x,s),
\end{array} \]
Dividing by $i \bar h(x,s) h(x,s)$, we will obtain that the last formula is
\[ \begin{array}{lll} 
2 \Theta_h'(x,s) =  \dfrac{1}{\bar h(x,s) h(x,s)}
\bbmatrix{1&\bar s} S^*(\lambda,x) \sigma_1  S(\lambda,x) \bbmatrix{1\\s} - (s + \bar s) = \\
= \dfrac{\lambda+\bar\lambda}{\bar h(x,s) h(x,s)} K_s(x,s) - (s + \bar s).
\end{array} \]
Finally, using part 2 we obtain that $\dfrac{|h(x,s)|^2 - \dfrac{\partial}{\partial x} K_S(s,x)}{i(s-\bar s)}=K_S(x)$,
which plugged into the last formula gives
\[ 2\Theta'_h(x) = \dfrac{\lambda+\bar\lambda}{\bar h(x,s) h(x,s)} \dfrac{|h(x,s)|^2 - \dfrac{\partial}{\partial x} K_S(s,x)}{i(s-\bar s)} - (s + \bar s) =
-\dfrac{s+\bar s}{|h(x,s)|^2} \dfrac{\partial}{\partial x} K_S(x,s).
\]

\qed

Important corollary of this theorem is a sort of independence of the formulas for $h(x,s)$ 
on the realization we choose for the vessel:
\begin{cor} Let $\mathfrak{K_V}$ and $\widetilde {\mathfrak{K_V}}$ be two vessels realizing the same potential 
$q(x)$. Suppose that $h(x,s) = |h(x,s)| e^{i\Theta_h(x,s)}$ and 
$\widetilde h(x,s) = |\widetilde h(x,s)| e^{i\widetilde \Theta_h(x,s)}$ are the corresponding functions, 
defined by \eqref{eq:h}. Then there exists a function $H(s)$, x-independent so that
\[ |h(x,s)|^2 = |\widetilde h(x,s)| H(s), \quad \Theta'_h(x,s) = \widetilde \Theta'_h(x,s).
\]
\end{cor}
\noindent\textbf{Proof:} From Theorem \ref{thm:UniquenessOfS} it follows that there exist two analytic
functions $a(\lambda), c(\lambda)$, which are zero at infinity such that for \eqref{eq:InCommPhi}
\[
 Y(\lambda) = I + \bbmatrix{a(\lambda)& \dfrac{i c(\lambda)}{\lambda} \\ c(\lambda) & a(\lambda)}
\]
it holds $\widetilde S(\lambda,x) = S(\lambda,x) Y(\lambda)$, for the corresponding transfer functions.
Then we calculate
\[ \begin{array}{lll}
\widetilde K_S(s,x)
= \dfrac{\bbmatrix{1 & \bar s}\widetilde S^*(\lambda,x) \sigma_1 \widetilde S(\lambda,x)\bbmatrix{1 \\ s} }{\lambda+\bar\lambda} = 
\dfrac{\bbmatrix{1 & \bar s} Y^*(\lambda) S^*(\lambda,x) \sigma_1 S(\lambda,x) Y(\lambda)\bbmatrix{1 \\ s} }{\lambda+\bar\lambda}.
\end{array} \]
Notice that
\[ Y(\lambda)\bbmatrix{1 \\ s} =  \bbmatrix{1+a(\lambda)& \dfrac{c(\lambda)}{s^2} \\ c(\lambda) & 1+a(\lambda)}
\bbmatrix{1 \\ s} = (1+a(\lambda) + \dfrac{c(\lambda)}{s}) \bbmatrix{1 \\ s},
\]
from where the result follows denoting $H(s) = |1+a(\lambda) + \dfrac{c(\lambda)}{s}|^2$ and the formula
\[ \widetilde K_S(s,x) = H(s) K_S(s,x).
\]
\qed

As a corollary of this theorem, we obtain a necessary conditions on $K_S(x,s)$ following from the existence
of the Jost solution $f(x,s)$:
\begin{cor} Suppose that we are given a vessel $\mathfrak {K_V}$, existing on a half-line
$\mathrm I=[x_0,\infty)$. Suppose also that for some value $s$, the Jost solution $f(x,s)$ exists, i.e. satisfies
condition \eqref{eq:fCond}. Then
\begin{enumerate}
	\item $\lim\limits_{x\rightarrow\infty} K_S(x,s)$ exists,
	\item $\lim\limits_{x\rightarrow\infty} \dfrac{\partial}{\partial x} K_S(x,s)$ exists,
\end{enumerate}
\end{cor}
\noindent\textbf{Proof:} From the Jost condition \eqref{eq:fCond} it follows that there exist two limits:
$\lim\limits_{x\rightarrow\infty} |h(x,s)|$ and $\lim\limits_{x\rightarrow\infty} \Theta_h(x,s)$. Then integrating
part 3 of Theorem \ref{thm:hisOK} it follows that
\[ 2(\lim\limits_{x\rightarrow\infty} \Theta_h(x,s) - \Theta_h(x_0,s)) = -\int\limits_{x_0}^\infty \dfrac{s+\bar s}{|h(y,s)|^2} \dfrac{\partial}{\partial y} K_S(y,s) dy.
\]
Thus the improper integral on the right hand side exists. As a result its integrand satisfies the necessary condition
of the convergence
\[ \lim\limits_{x\rightarrow\infty} \dfrac{s+\bar s}{|h(x,s)|^2} \dfrac{\partial}{\partial x} K_S(x,s) = 0.
\]
Dividing next part 2 of Theorem \ref{thm:hisOK} we find that
\[ 1 = \dfrac{\dfrac{\partial}{\partial x} K_S(x,s)}{|h(x,s)|^2} + i(s-\bar s)\dfrac{K_S(x,s)}{|h(x,s)|^2}
\]
Taking $x$ approaching infinity on both sides we find that
\[ \lim\limits_{x\rightarrow\infty} i(s-\bar s) \dfrac{K_S(x,s)}{|h(x,s)|^2} = 1,
\]
or that (since $h(\infty,s) \neq 0)$
\[ \lim\limits_{x\rightarrow\infty} K_S(x,s) = \dfrac{|h(\infty,s)|^2}{i(s-\bar s)}.
\]
Finally, using again part 2 of Theorem \ref{thm:hisOK}
\[ \dfrac{\partial}{\partial x} K_S(s,x) = |h(x,s)|^2 -  i(s-\bar s) K_S(x),
\]
we obtain that since the right hand side has a limits as $x$ approaches infinity, so does the left.
\qed

\section{Applications}
\subsection{\label{sec:SLDissImagine}SL vessels with a spectrum on the imaginary positive axis.}
When the spectrum of the operator $A$ is on the imaginary axis, it means that $m(A)$, defined in \eqref{eq:mA} is zero:
\[ m(A) = \max \{ \Im \lambda \mid \lambda\in\SPEC(A)\} = 0.
\]
Let us consider a vessel (see Definition \ref{def:KV})
\[\mathfrak{K_V} = (A, B(x), \mathbb X(x); \sigma_1(x), 
\sigma_2(x), \gamma(x), \gamma_*(x);
\mathcal{H}, \mathcal{E};\mathrm I), 
\]
which has an additional restriction, identical to the classical case,
namely, the operator $A$ has all its spectrum on the imaginary positive axis: 
\begin{equation} \label{eq:Aclassic}
\SPEC A \subseteq i\mathbb R^+.
\end{equation}
We can find more accurate bounds then in Lemma \ref{lemma:PosSig2} on norms of vessel operators. Starting from 
\[ \begin{array}{lll}
B(x) = \dfrac{1}{2\pi i} \oint (\lambda I - A)^{-1} B_0 \Phi^t (\lambda,x,x_0) d\lambda
\end{array} \]
and denoting $K = \sqrt{-iA}>0$ (exists from \eqref{eq:Aclassic}), we obtain 
\[ \begin{array}{llll}
B(x) & = \cos(K x) B(x_0) + i K\sin(K x) B(x_0) \bbmatrix{0&1\\0&0} +
 (K)^{-1}\sin(K x) B(x_0) \bbmatrix{0&0\\1&0}, \\
 \mathbb X(x) & = \mathbb X(x_0) + \int_{x_0}^x \{ \cos(K y) B_0 \sigma_2 B_0^* \cos(K y) +
K^{-1} \sin(K y) B_0 \bbmatrix{0&0\\0&1} B_0^* \sin(K y) K^{-1} - \\
&  - i \cos(K y) B_0  \bbmatrix{0&1\\0&0} B_0^*  \sin(K y) K^{-1} +
 i K^{-1} \sin(K y) B_0  \bbmatrix{0&0\\1&0} B^*(x_0) \cos(K y) \} dy.
\end{array} \]
From where we obtain the following bounds:
\begin{thm} \label{thm:QBoundGen}
Let $\mathfrak {K_V}$ be a vessel, existing on $[x_0,\infty)$, for which $A$ satisfies condition \eqref{eq:Aclassic}. 
Then on the interval $[x_0,\infty)$ the following bounds hold
\begin{enumerate}
	\item $\| B(x) \| \leq B_1$ and $\|\mathbb X(x)\| \leq \|\mathbb X(x_0)\| + B_2(x-x_0)$
		for some $B_1, B_2$,
	\item $\TR (\mathbb X(x) - \mathbb X(x_0)) \geq T_1 x + T_2$ for some constants $T_1>0, T_2$.
\end{enumerate}
\end{thm}
\noindent\textbf{Proof: 1.} Examining the formulas above, 
$B(x)$ is determined using formulas of the form $\cos(Ky), \sin(Ky)$, which are
bounded by $1$ for a positive $K$. Integrating further these expressions, we will obtain the bound for 
$\|\mathbb X(x)\|$.

\noindent\textbf{2.} Consider the following calculation
\[ \begin{array}{lllll} 
\TR (\mathbb X(x) - \mathbb X_0) = \int_{x_0}^x \{ \cos(K y) B_0 \sigma_2 B_0^* \cos(K y) +
K^{-1} \sin(K y) B_0 \bbmatrix{0&0\\0&1} B_0^* \sin(K y) K^{-1} - \\
  - i \cos(K y) B_0  \bbmatrix{0&1\\0&0} B_0^*  \sin(K y) K^{-1} +
 i K^{-1} \sin(K y) B_0  \bbmatrix{0&0\\1&0} B^*(x_0) \cos(K y) \} dy ) = \\
= \int_{x_0}^x \TR (B_0 \sigma_2 B_0^* \cos(K y)  \cos(K y)+
B_0 \bbmatrix{0&0\\0&1} B_0^* \sin(K y) K^{-1}K^{-1} \sin(K y) - \\
- i B_0  \bbmatrix{0&1\\0&0} B_0^* \sin(K y) K^{-1}\cos(K y) +
i  B_0  \bbmatrix{0&0\\1&0} B^*(x_0) \cos(K y)K^{-1} \sin(K y) ) dy ) = \\
= \int_{x_0}^x \TR (B_0 \sigma_2 B_0^* \dfrac{1 + \cos(2K y)}{2}+
B_0 \bbmatrix{0&0\\0&1} B_0^* \dfrac{1-\cos(2K y)}{2} K^{-2} - \\
- i B_0  \bbmatrix{0&1\\0&0} B_0^* \sin(K y) K^{-1}\cos(K y) +
i  B_0  \bbmatrix{0&0\\1&0} B^*(x_0) \cos(K y)K^{-1} \sin(K y) ) dy ) = \\
= \int_{x_0}^x \TR (B_0 \sigma_2 B_0^* + B_0 \bbmatrix{0&0\\0&1} B_0^* K^{-2}) dy + \\
\quad + \int_{x_0}^x \TR(B_0 \sigma_2 B_0^*\dfrac{\cos(2K y)}{2}+
B_0 \bbmatrix{0&0\\0&1} B_0^* \dfrac{-\cos(2K y)}{2} K^{-2} - \\
- i B_0  \bbmatrix{0&1\\0&0} B_0^* \sin(K y) K^{-1}\cos(K y) +
i  B_0  \bbmatrix{0&0\\1&0} B^*(x_0) \cos(K y)K^{-1} \sin(K y) ) dy ) = \\
= T_1 x + T_2(x),
\end{array} \]
where $T_1>0$ and $T_2(x)$ is uniformly bounded for all $x$. So defining $T_2=-\inf T_2(x)$, we will obtain the
desired result.
\qed

More accurate bounds are obtained in the dissipative case and are presented in the next theorem.
\begin{thm} \label{thm:QBoundDiss}
Let $\mathfrak {K_V}$ be a dissipative vessel, existing on $[x_0,\infty)$, for which $A$ satisfies condition \eqref{eq:Aclassic}. 
Then, on the interval $[x_0,\infty)$ the following bounds hold
\begin{enumerate}
	\item $\tau(x) \geq T_1 (x-x_0) + T_2$ for some constants $T_1, T_2 > 0$,
	\item $\|\mathbb X^{-1}(x)\| \leq \dfrac{1}{\tau(x)} \leq \dfrac{1}{T_1 (x-x_0) + T_2}$,
	\item $|q(x)| \leq \dfrac{Q}{(x-x_0)}$ for some $Q$ and for big enough $x$.
\end{enumerate}
\end{thm}
\noindent\textbf{Proof: 1.}  Notice that using Theorem \ref{thm:X0isI} we may suppose that $\mathbb X_0=I$.
Then, using the formula
\[ \det (I + F) > 1 + \TR F
\]
for a trace class positive operator $F$ and part 2 of Theorem \ref{thm:QBoundGen},
we will obtain the bound for $\tau(x)$.

\noindent\textbf{2.} Follows from the following chain of inequalities
\[ \begin{array}{llll}
\| \mathbb X^{-1}(x) \| &\leq \| I \| + \| \mathbb X^{-1}(x) - I \| \leq
1 + \| \mathbb X^{-1}(x) - I \|_1 = 
1 + \TR (\mathbb X^{-1}(x) - I) \leq \\
& \leq \det(\mathbb X^{-1}(x)) = \dfrac{1}{\det \mathbb X(x)} = \dfrac{1}{\tau(x)} \leq \dfrac{1}{T_1 (x-x_0) + T_2}.
\end{array} \]

\noindent\textbf{3.} Recall the formula \eqref{eq:qbeta} for the potential.
Differentiating it and using vessel equations we will arrive to
\[ \begin{array}{lll}
-\dfrac{1}{2} q(x)  & = \dfrac{d}{dx}\TR(\mathbb X'(x) \mathbb X^{-1}(x)) = \\
& = \TR(\bbmatrix{0&i\\-i&0} B^*(x)\mathbb X^{-1}(x) B(x)) -
(\bbmatrix{1&0} B^*(x)\mathbb X^{-1}(x) B(x)\bbmatrix{1\\0})^2.
\end{array}\]
using part 1 of Theorem \ref{thm:QBoundGen} and part 2 of this Theorem we obtain that 
\[ |q(x)| \leq Q_1 \| \mathbb X^{-1}\|_\infty + Q_2 \| \mathbb X^{-1}\|_\infty^2 \leq
\dfrac{Q_1}{T_1 (x-x_0) + T_2} + \dfrac{Q_2^2}{(T_1 (x-x_0) + T_2)^2} \leq \dfrac{Q}{(x-x_0)}
\]
for big enough $x$.
\qed
\subsection{\label{sec:ArbCurve}Models of vessels with a spectrum on a symmetric 
curve $\Gamma$.}
Let us construct a vessel, for which the operator $A$ is diagonal.
\begin{enumerate}
\item Suppose that we are given a smooth curve in the complex domain, parametrized by
	\[ \Gamma = \{ \mu(t) \mid a\leq t \leq b \}
	\]
	Let us also suppose that the curve is symmetric with respect to the imaginary axis, i.e.
	\linebreak 	$\Gamma = - \Gamma^*$. The reason why we need to require it is the symmetry
	condition \eqref{eq:Symmetry}, which means that if $\lambda_0$ is an analytic point of 
	$S(\lambda,x)$ and $S(\lambda,x) = \sum \alpha_n(x) (\lambda-\lambda_0)^n$, then
	\[ S^*(-\bar\lambda,x) = \sum \sigma_1(x) \alpha_n(x) \sigma_1^{-1}(x) (\lambda-\lambda_0)^n,
	\]
	and consequently, taking the conjugate
	\[ S(-\bar\lambda,x) = 
	\sum \sigma_1^{-1}(x) \alpha_n^*(x) \sigma_1(x) (-1)^n (-\bar\lambda + \bar \lambda_0)^n,
	\]
	where substituting $-\bar\lambda=\mu$, we will obtain that $S(\mu,x)$ is analytic
	at $-\bar\lambda_0$. Particularly, it can't be a singular point of $S(\lambda,x)$. 
	Since we will verify that the singularities occur on the
	curve $\Gamma$, this means that the curve must be symmetric.
	The inner space is defined as
	\[ \mathcal H = L^2(\Gamma)=\{ f(\mu) \mid \int\limits_a^b |f(\mu(t))|^2 dt < \infty \}.
	\]
	with the inner product
	\[ \langle f(\mu), g(\mu) \rangle_\mathcal H = \int_a^b g^*(\mu(t)) f(\mu(t)) dt.
	\]
\item Define the operator $A$ as the multiplication operator on a function $\mu$:
	\[ A f(\mu) = \mu f(\mu)
	\]
\item Then $B(x)$ is a solution of \eqref{eq:DB}
	\[ 0  =  \frac{d}{dx} (B(x)\sigma_1(x)) + A B(x) \sigma_2(x) + B(x) \gamma(x).
	\]
	For example, in the SL case, we shall obtain that $B(x)$ is an operator by multiplication
	on $B(\mu,x) = \bbmatrix{c_1(\mu) & c_2(\mu)} \Phi(\mu,x,x_0)$,	for "good" (satisfying
	\eqref{eq:BIsSchwz})	functions $c_1(\mu), c_2(\mu)$. In order to verify that the
	singularities of the final transfer function will occur exactly on the curve, we need the
	minimality condition, which must be satisfied at $x_0$ by Theorem \ref{thm:PermMin}. Since
	$A = \mu$, it is enough to demand that $B(\mu,x_0) \neq 0$ for all $\mu$, since then the
	condition \eqref{eq:minCond} is immediate.
	
	Notice that from the definition it follows that the adjoint of $B(x)$ is
	\begin{equation} \label{eq:DiagB*}
	B^*(x) f(\mu) = \int\limits_a^b B^*(\mu(t),x) f(\mu(t)) dt.
	\end{equation}
	\item Define the operator $\mathbb X(x)$ as follows
	\begin{equation} \label{eq:DiagX}
	\mathbb X(x) f(\mu) = \int\limits_a^b \dfrac{B(\mu,x) \sigma_1 B^*(\delta(t),x)}{\mu+\bar\delta(t)} dt.
	\end{equation}
\end{enumerate}
Notice that in order to obtain a well-defined operator, we have to verify that the integral converges. For this to hold,
we need to verify that for values of $\delta$, where $\mu + \bar\delta(t) = 0$, it holds that 
\[
\dfrac{B(\mu,x) \sigma_1 B^*(\delta(t),x)}{\mu+\bar\delta(t)}
\]
is integrable. One can demand for that that the H\"older condition \cite{bib:Mushel} is satisfied, for example. Moreover, if the curve $\Gamma$ is unbounded, we have to choose
$B_0 =B(a) = B(\mu,a)$ such that $\|\mu^n B(\mu)\| \leq C^n$ for a constant $C$. For this it is enough to choose $B(\mu)$ to be a Schwartz function on the curve $\Gamma$.

\begin{thm} The following collection
\begin{equation} \label{eq:DiagV}
\mathfrak{V}_d = (A = \mu, B(\mu,x), \mathbb X(x); \sigma_1, \sigma_2, \gamma, \gamma_*(x);
\mathcal{H}=L^2(\Gamma),\mathcal E), 
\end{equation}
for $\gamma_*(x)$, defined by the linkage condition \eqref{eq:Linkage}, is a vessel.
\end{thm}
\noindent\textbf{Proof:} Equation \eqref{eq:DB} is satisfied by the construction of $B(\mu,x)$. Lyapunov equation 
\eqref{eq:XLyapunov} follows from the definitions
\begin{multline*}
(A_1 \mathbb X(x) + \mathbb X(x) A^*) f(\mu) = \\
= \mu \int\limits_a^b \dfrac{B(\mu,x) \sigma_1 B^*(\delta(t),x)}{\mu+\bar\delta(t)} f(\delta(t)) dt +
\int\limits_a^b \dfrac{B(\mu,x) \sigma_1 B^*(\delta(t),x)}{\mu+\bar\delta(t)} a^*(\delta(t)) f(\delta(t)) dt = \\
= \int\limits_a^b \dfrac{B(\mu,x) \sigma_1 B^*(\delta(t),x)}{\mu+\bar\delta(t)} (\mu+\bar\delta(t)) f(\delta(t)) dt = \\
= \int\limits_a^b B(\mu,x) \sigma_1 B^*(\delta(t),x) f(\delta(t)) dt = B(x) \sigma_1 B^*(x) f(\mu).
\end{multline*}
Similarly, equation \eqref{eq:DX} holds. Finally, \eqref{eq:Linkage} serves to define $\gamma_*(x)$.
\qed

Finally notice that the transfer function of this vessel is
\[ \begin{array}{ll}
S(\lambda,x) & = I - B^*(x) \mathbb X^{-1}(x) (\lambda I - A)^{-1} B(x) \sigma_1(x) \\
& = I - \int\limits_a^b \dfrac{B^*(\mu(t),x) \mathbb X^{-1}(x) B(\mu(t),x)}{\lambda - \mu(t)} \sigma_1(x) dt
\end{array} \]
and has singularities (jumps) exactly on the curve $\Gamma$ by the minimality condition.

\subsection{\label{sec:NLS}NLS equations}
The first part of the article suggests that one can construct vessels with a prescribed spectrum for a wider class
of vessel parameters. For example one can study Non Linear Schr\" odinger (NLS) equations presented
by A.P. Fordy, P.P Kulish in \cite{bib:NLS}. The classical NLS equation corresponds to the following parameters:
\begin{defn} The Non-Linear Shr\" odinger equation parameters are defined to be
\[ \begin{array}{lll}
 \sigma_1 = \bbmatrix{1&0\\0&1},\quad
\sigma_2 = \dfrac{1}{2} \bbmatrix{1&0\\0&-1},
\quad \gamma(x) =\bbmatrix{0 & 0 \\0 & 0}, \\
\gamma_*(x) = \bbmatrix{0&\beta(x) \\-\beta^*(x)& 0}.
\end{array} \]
\end{defn}
\noindent and the regularity assumptions can be taken as in the SL case.
The output compatibility conditions take the form of the classical non linear Schr\" odinger 
equation with the spectral parameter $i\lambda$
\[ \dfrac{\partial}{\partial x} u(x,\lambda) = (i \lambda A + Q(x)) u(x,\lambda),
\]
where 
\[ I=\sigma_1, \quad A=-\dfrac{1}{2}\bbmatrix{i&0\\0&-i} = -i \sigma_2, \quad Q(x) = -\gamma_*(x).
\]
A more complicated example is \cite[3.19]{bib:NLS} as follows
\[ \dfrac{\partial}{\partial x} \bbmatrix{\phi_1\\\phi_2\\\phi_3\\\phi_4} =
\bbmatrix{
\dfrac{1}{2} i\lambda & 0 & q_1 & q_2\\
0 & \dfrac{1}{2} i\lambda & q_4 & q_3\\
-q_1^* & -q_4^* & -\dfrac{1}{2} i\lambda & 0\\
-q_2^* & -q_3^* & 0 & -\dfrac{1}{2} i\lambda \\
} \bbmatrix{\phi_1\\\phi_2\\\phi_3\\\phi_4}
\]
and corresponds to the following vessel parameters (changing the spectral parameter $\lambda$ to $-i\lambda$.
\[ \sigma_1 = I, \quad \sigma_2 = \operatorname{diag} \{ \dfrac{1}{2}, \dfrac{1}{2}, -\dfrac{1}{2},-\dfrac{1}{2}\}, \quad \gamma=0,
\quad \gamma_*(x) =
\bbmatrix{
0 & 0 & q_1 & q_2\\
0 & 0 & q_4 & q_3\\
-q_1^* & -q_4^* & 0 & 0\\
-q_2^* & -q_3^* & 0 & 0 \\
}.
\]
One can use models, presented in section \ref{sec:ArbCurve} for the study of the Scattering Theory of these
NLS equations.

It is important to notice that the connection between the corresponding to NLS parameters vessels and
the Lie algebras appearing in \cite{bib:NLS} is of great interest.

\subsection{\label{sec:CanSys}Canonical Systems}
Starting from a standard model of a canonical system \cite{bib:FaddeyevII}
\[ [J \dfrac{d}{dx} + Q(x)] \phi(x,k) = k \phi(x,k)
\]
where $J = \bbmatrix{0&1\\-1&0}$, $Q(x) = \bbmatrix{p(x)&q(x)\\q(x)&-p(x)}$ notice that
multiplying this equation by $i$, we will obtain a differential equation,
which fits the setting of a vessel:
\[ [\bbmatrix{0&i\\-i&0} \dfrac{d}{dx} - ik \phi(x,k) - 
\bbmatrix{-ip(x)&-iq(x)\\-iq(x)& ip(x)}] \phi(x,k) = 0.
\]
\begin{defn} The canonical system vessel parameters are
\[ \sigma_1 = \bbmatrix{0&i\\-i&0}, \sigma_2 = I, \gamma=0,
\gamma_*(x)=\bbmatrix{-ip(x)&-iq(x)\\-iq(x)& ip(x)}.
\]
\end{defn}
\noindent Requiring \cite[2.2]{bib:FaddeyevII}
\[ \int\limits_{-\infty}^\infty |q(x)| dx < \infty, 
\int\limits_{-\infty}^\infty |p(x)| dx < \infty
\]
we will obtain a vessel with a spectrum on a cut of the imaginary positive axis,
imitating the construction for SL case and using the formulas from
\cite[section 2]{bib:FaddeyevII}. The general case will produce an interesting class of
potentials $Q(x)$ in this case too.

\noindent\textbf{Acknowledgment}
I would like to thank Ronald K. Perlin from Drexel university for introducing and explaining me NLS equations.


\bibliographystyle{alpha}
\bibliography{biblio}

\begin{thebibliography}{{M}us41}

\bibitem[AD]{bib:KreinReal}
H.S.V. de~{S}noo A.~{D}iksma, H.~{L}anger.
\newblock Representations of holomrphic operator functions by means of
  resolvents of unitary or self-adjoint operators in {K}rein spaces.
\newblock {\em Operator Theory: Adv. and App.}, 24:123--143.
\newblock Birkhauser Verlag, Berlin.

\bibitem[AF83]{bib:NLS}
P.P.~{K}ulish A.P.~{F}ordy.
\newblock Non linear schr{\" o}dinger equations and simple lie algebras.
\newblock {\em Communications in Math. Phys.}, 89:427--443, 1983.

\bibitem[AM09]{bib:amv}
D.~{A}lpay A.~Melnikov, V.~Vinnikov.
\newblock Un algorithme de schur pour les fonctions de transfert des syst\`emes
  surd\'etermin\'es invariants dans une direction.
\newblock {\em Comptes-Rendus math\'ematiques (Paris)}, 347(13--14):729--733,
  2009.

\bibitem[{A}ro50]{bib:AronsjKErnels}
N.~{A}ronszajn.
\newblock Theory of reproducing kernels.
\newblock {\em Tran. of AMS}, 68(3):337–--404, 1950.

\bibitem[BGR90]{bib:bgr}
J.~{B}all, I.~{G}ohberg, and L.~{R}odman.
\newblock {\em Interpolation of rational matrix functions}.
\newblock Operator Theory: Advances and Applications. Birkh{\" a}user Verlag,
  Basel, 1990.

\bibitem[Bi71]{bib:Brodskii}
M.S. {B}rodski\u i.
\newblock {\em Triangular and Jordan representations of linear operators}.
\newblock translations of AMS, 1971.

\bibitem[CL55]{bib:CoddLev}
E.A {C}oddington and N.~{L}evinson.
\newblock {\em Theory of ordinary differential equation}.
\newblock Mc-Graw Hill, 1955.

\bibitem[{C}ru55]{bib:Crum}
M.M. {C}rum.
\newblock Associated {S}turm-{L}iouville systems.
\newblock {\em Quart. J. Math. Oxford Ser.}, 6(2):121--127, 1955.

\bibitem[{F}ad63]{bib:FadeevInv}
L.D. {F}adeev.
\newblock The inverse problem in the quantum theory of scattering.
\newblock {\em Journal of Mathematical Physics}, 4(1):72--104, 1963.

\bibitem[{F}ad74]{bib:FaddeyevII}
L.D. {F}adeev.
\newblock The inverse problem in the quantum theory of scattering, {II}.
\newblock {\em Itogi Nauk. i Techn.}, 4:93--180, 1974.

\bibitem[Fuc07]{bib:Fuchs}
R.~Fuchs.
\newblock {\" U}ber lineare homogene differentialgleichungen zweiter ordnung
  mit drei im endlichen gelegenen wesentlich singul{\" a}ren stellen
  ({G}erman).
\newblock {\em Math. Ann.}, 63(3):301--321, 1907.

\bibitem[Gar12]{bib:Garnier}
R.~Garnier.
\newblock Sur des {\' e}quations diff{\' e}rentielles du troisi{\` e}me ordre
  dont l'int{\' e}grale g{\' e}n{\' e}rale est uniforme et sur une classe d'{\'
  e}quations nouvelles d'ordre sup{\' e}rieur dont l'int{\' e}grale g{\' e}n{\'
  e}rale a ses points critiques fixes ({F}rench).
\newblock {\em Ann. Sci. {\' E}cole Norm. Sup.}, 29(3):1--126, 1912.

\bibitem[{H}el74]{bib:helton}
J.W. {H}elton.
\newblock Discrete time systems, operator models, and scattering theory.
\newblock {\em Journal of Functional analysis}, 16:15--38, 1974.

\bibitem[IG69]{bib:GKintro}
M.~{K}rein I.~{G}ohberg.
\newblock {\em Introduction to the theory of linear non-selfadjoint operators}.
\newblock translations of AMS, 1969.

\bibitem[IMG51]{bib:GL}
B.~M.~{L}evitan I.~M.~{G}elfand.
\newblock On the determination of a differential equation from its spectral
  function ({R}ussian).
\newblock {\em Izvestiya Akad. Nauk SSSR. Ser. Mat.}, 15, 1951.

\bibitem[{J}os47]{bib:Jost}
Res {J}ost.
\newblock Bemerkungen zur matematischen theorie der z{\" a}hler ({G}erman).
\newblock {\em Helvetica Phys. Acta}, 20:173--–182, 1947.

\bibitem[KR69]{bib:Kalman}
{A}rbib~M.A. {K}alman R.E., {F}alb~P.L.
\newblock {\em Topics in mathematical system theory}.
\newblock McGraw-Hill, 1969.

\bibitem[{K}re55]{bib:Krein1}
M.~G. {K}rein.
\newblock On determination of the potential of a particle from its $s$-function
  ({R}ussian).
\newblock {\em Dokl. Akad. Nauk SSSR}, 105:433–--436, 1955.

\bibitem[{L}ev49]{bib:Levinson}
N.~{L}evinson.
\newblock On the uniqueness of the potential in a {S}chr{\" o}dinger equation
  for a given asymptotic phase.
\newblock {\em Danske Vid. Selsk. Mat.-Fys. Medd.}, 25(9):29--, 1949.

\bibitem[{L}io95]{bib:Liouville}
R.~{L}iouville.
\newblock Sur les \'equations de la dynamique ({F}rench).
\newblock {\em Acta Math.}, 19(1):251--283, 1895.

\bibitem[LP67]{bib:LaxPhil}
P.D. {L}ax and R.S. {P}hilips.
\newblock {\em Scattering theory}.
\newblock Academic Press, 1967.

\bibitem[Ls78]{bib:defVess}
M.S. {L}iv\v sic.
\newblock Commuting nonselfadjoint operators and solutions of systems of
  partial differential equations generated by them, ({R}ussian).
\newblock {\em Soobshch. Akad. Nauk Gruzin. SSSR}, 91(2):281--284, 1978.

\bibitem[Ls01]{bib:Vortices}
M.S. {L}iv\v sic.
\newblock Vortices of 2d systems.
\newblock {\em Operator Theory: Advances and Applications}, 123:7--41, 2001.

\bibitem[MB58]{bib:BL}
M.S.~{L}iv{\v s}ic M.S.~{B}rodskii.
\newblock Spectral analysis of non-self-adjoint operators and intermediate
  systems ({R}ussian).
\newblock {\em Uspehi Mat. Nauk (N.S.)}, 13(1 (79)):3--85, 1958.

\bibitem[Mc50]{bib:Marchenko}
V.~A. {M}ar\v cenko.
\newblock Concerning the theory of a differential operator of the second order
  ({R}ussian).
\newblock {\em Doklady Akad. Nauk SSSR. (N.S.)}, 72:457--–460, 1950.

\bibitem[Mc77]{bib:MarchenkoSL}
V.A. {M}ar\v cenko.
\newblock {\em Sturm Liouville operators and their applications}.
\newblock Naukova Dumka, Kiev, 1977.

\bibitem[{M}el]{bib:SLVessels}
A.~{M}elnikov.
\newblock Finite dimensional sturm liouville vessels and their tau functions.
\newblock accepted to IEOT.

\bibitem[{M}el09]{bib:MyThesis}
A.~{M}elnikov.
\newblock {\em Overdetermied $2D$ systems invariant in one direction and their
  transfer functions}.
\newblock PhD thesis, Ben Gurion University, 2009.

\bibitem[{M}us41]{bib:Mushel}
N.I. {M}uskhelishvili.
\newblock Application of integrals of cauchy type to a class of singular
  integral equations.
\newblock {\em Trudy Tbilissi. Mat. Obschestva}, 10:1--44, 1941.

\bibitem[MVa]{bib:MelVin1}
A.~{Melnikov} and V.~{Vinnikov}.
\newblock Overdetermined 2{D} systems invariant in one direction and their
  transfer functions.
\newblock http://arXiv.org/abs/0812.3779.

\bibitem[MVb]{bib:MelVinC}
A.~{Melnikov} and V.~{Vinnikov}.
\newblock Overdetermined conservative 2{D} systems, invariant in one direction
  and a generalization of {P}otapov's theorem.
\newblock http://arxiv.org/abs/0812.3970.

\bibitem[MVc]{bib:SchurVessels}
D.~{Alpay}~A. {Melnikov} and V.~{Vinnikov}.
\newblock On the class $\mathbf {SI}$ of conservative functions intertwining
  solutions of linear differential equations.
\newblock http://arxiv.org/abs/0912.2014.

\bibitem[MV60]{bib:AgMarch}
{A}granovich~Z.S. {M}archenko V.A.
\newblock {\em Inverse problem of {S}cattering theory ({R}ussian)}.
\newblock KGU, Kharkov, 1960.

\bibitem[ND88]{bib:DanSchw}
J.~T.~{S}chwartz N.~{D}unford.
\newblock {\em Linear Operators, General Theory}.
\newblock Wiley-Interscience, 1988.

\bibitem[{P}ot55]{bib:Potapov}
V.P. {P}otapov.
\newblock On the multiplicative structure of $j$-nonexpanding matrix functions
  (russian).
\newblock {\em Trudy Moskov. Mat. Obschestva}, 4:125--236, 1955.
\newblock English transl.: AMS Translations (2), 15 (1960), 131--243.

\bibitem[{P}ov50]{bib:Povzner}
A.~{P}ovzner.
\newblock On differential equations of {S}turm-{L}iouville type on a half-axis.
\newblock {\em Amer. Math. Soc. Translation}, 79(5), 1950.

\bibitem[{S}ch08]{bib:Schlez}
L.~{S}chlesinger.
\newblock Sur la solution du probl{\'e}me de {R}iemann ({F}rench), 1908.

\bibitem[{S}tu36]{bib:Sturm}
C.~{S}turm.
\newblock Sur les {\' e}quations diff{\' e}rentielles lin{\' e}aires du second
  ordre ({F}rench), 1836.

\bibitem[{Z}el98]{bib:Zelikin}
M.I. {Z}elikin.
\newblock {\em Unimodular spaces and Riccati equation in variational calculus}.
\newblock Factorial, 1998.

\end{thebibliography}

\end{document}